\documentclass[11pt,a4paper]{article}

\usepackage{iftex}
\ifPDFTeX
  \usepackage[T1]{fontenc}
  \usepackage[utf8]{inputenc}
  \usepackage{lmodern}
\else
  \usepackage{fontspec}
\fi
\usepackage{amsmath,amssymb,amsthm,mathtools}
\usepackage[a4paper,margin=2.6cm]{geometry}
\usepackage{enumitem}
\usepackage{booktabs}
\usepackage[hidelinks]{hyperref}

\theoremstyle{plain}

\newtheorem{theorem}{Theorem}[section]
\newtheorem{proposition}[theorem]{Proposition}
\newtheorem{lemma}[theorem]{Lemma}
\newtheorem{corollary}[theorem]{Corollary}
\theoremstyle{definition}
\newtheorem{definition}[theorem]{Definition}
\newtheorem{example}[theorem]{Example}
\newtheorem{remark}[theorem]{Remark}
\newtheorem{question}[theorem]{Question}

\numberwithin{equation}{section}
\setlist[enumerate]{leftmargin=2.2em}

\newcommand{\Z}{\mathbb{Z}}

\newcommand{\F}{\mathbb{F}}
\newcommand{\Mat}{\mathbb{M}}
\newcommand{\Tri}{\mathbb{T}}
\newcommand{\dd}{\mathbf{d}}
\newcommand{\rr}{\mathbf{r}}
\newcommand{\gl}{\mathfrak{gl}}
\newcommand{\cG}{\mathcal{G}}
\newcommand{\cD}{\mathcal{D}}
\newcommand{\cJ}{\mathcal{J}}
\newcommand{\Nil}{\operatorname{Nil}_*}

\DeclareMathOperator{\chr}{char}
\DeclareMathOperator{\spn}{span}
\DeclareMathOperator{\dl}{dl}
\DeclareMathOperator{\jl}{jl}
\newcommand{\dlL}{\dl^{L}}

\begin{document}

\begin{center}
{\LARGE\bfseries Jordan products and Lie solvability\\[4pt] of matrix rings\par}

\vspace{14pt}
{\large Mai Hoang Bien\textsuperscript{a,b}, Nguyen Thi Mong Tuyen\textsuperscript{c,$*$} and Le Anh Vu\textsuperscript{d,b}\par}

\vspace{6pt}
{\small\textsuperscript{$*$}Corresponding author\par}
\end{center}

\vspace{2pt}
\begin{abstract}
Let $N$ be an associative ring, not necessarily commutative or unital. We show that the matrix ring $\Mat_2(N)$ is Lie solvable if and only if the Jordan ideal series of $N$, built from the products $xy+yx$, terminates, and that the Lie derived length of $\Mat_2(N)$ lies between the length $\jl(N)$ of this series and $3\jl(N)$, the constant $3$ being best possible. When $2N=0$ the series is the strong Lie derived series of $N$, and when $2$ is invertible on $N$ its termination amounts to nilpotency, so that the characteristic hypothesis of the classical commutative case becomes an intrinsic condition on $N$. For $n\ge3$, by contrast, $\Mat_n(N)$ is Lie solvable exactly when $N$ is nilpotent, in every characteristic. The proof uses the three-term grading of $\gl_2(N)$, whose components of degree $\pm1$ carry the Jordan triple product of $N$. We also compute the derived series of $\gl_2(\Z)$ in closed form and show that, in a Lie solvable ring of derived length $d$, the idempotent $e_{11}$ of a system of $2\times2$ matrix units has additive order $2^m$ with $m\le2\lfloor2^{d-1}/3\rfloor$, a sharp bound. Applications concern inverse subsemigroups, group algebras and Steinberg algebras.

\end{abstract}

\section{Introduction}\label{s:intro}

Every associative ring $R$ is a Lie ring under the commutator $[x,y]=xy-yx$, and the interplay between the two structures has been studied since the work of Jennings \cite{Je42,Je47} and Herstein \cite{Her69}. For Lie solvable rings the fundamental structure theorem is due to Zalesski\u{\i} and Smirnov \cite{ZS82}, who proved that the ideal generated by $[[[R,R],[R,R]],R]$ is nilpotent (see \cite[Theorem~2.5]{SU17}, and \cite[Theorem~2.1]{SS85} for an explicit bound when $3$ is invertible). Passi, Passman and Sehgal \cite{PPS73} characterized the Lie solvable group algebras, whose Lie derived lengths were studied in \cite{CS07,Ju06,Sha92,Sp08}; for group algebras, Lie solvability and strong Lie solvability coincide in odd characteristic but not in characteristic $2$ \cite{Sp08}, a first sign of the special role of the prime $2$. Lie solvability has also been investigated for matrix algebras \cite{MSW12,WZ19}, tensor products \cite{Ri96}, finitary incidence algebras \cite{KS24} and Leavitt path algebras \cite{NZ22}. The Lie structure of Leavitt path algebras has attracted considerable attention: simple Lie algebras arising from Leavitt path algebras and Steinberg algebras were studied in \cite{AM12,Na20}, and commutator Leavitt path algebras in \cite{Me13}. Khanh and Danh \cite{KD25} determined when the Lie algebra of skew-symmetric elements and the Jordan algebra of symmetric elements of a Leavitt path algebra, with respect to its standard involution, are solvable; along the way they showed that, in characteristic $2$, the skew-symmetric $2\times2$ matrices over an integral domain with involution form a Lie solvable algebra of derived length at most $3$ \cite[Proposition~3.1]{KD25}. Throughout this theory one algebra plays an exceptional role: the algebra $\Mat_2(K)$ of $2\times2$ matrices over a field $K$ of characteristic $2$. It is the only full matrix algebra of size at least $2$ over a field that is Lie solvable \cite[Lemma~1.2]{PPS73}. Prime rings of characteristic $2$ that are $4$-dimensional over their extended centroid \cite{BMM96}, such as $\Mat_2(K)$, form the exceptional case in the Lie structure theory of prime rings \cite{Her69,LM72}. Moreover, for an arbitrary field $K$, the identity of Lie solvability of length $d$ holds in $\Mat_2(K)$ exactly when $\chr K=2$ and $d\ge3$; in all other cases it is a non-matrix identity, that is, one that fails in $\Mat_2(K)$. This paper is devoted to the following question, which asks how far this exception extends when the field is replaced by an arbitrary ring.

\begin{question}\label{q:main}
For which associative rings $N$, not necessarily commutative or unital, is the matrix ring $\Mat_2(N)$ Lie solvable, and what controls its Lie derived length?
\end{question}

Only matrices of size $2$ lead to a genuinely new question. For pairwise distinct indices $i,j,l$ one has $[xE_{ij},yE_{jl}]=xyE_{il}$, and for $n\ge3$ this identity alone forces all powers of $N$ into the derived series of $\Mat_n(N)$: whatever the characteristic, $\Mat_n(N)$ is Lie solvable exactly when $N$ is nilpotent, and its derived length is then determined by the index of nilpotency of $N$ (Proposition~\ref{prop:Mn}). For $n=2$ no third index is available, and the Lie structure of $\Mat_2(N)$ detects Jordan products instead, as explained below.

For a commutative algebra $S$ over a field, the answer is contained in a lemma of Passi, Passman and Sehgal \cite[Lemma~1.2]{PPS73}: if $S$ is not nilpotent, then $\Mat_2(S)$ is Lie solvable if and only if the characteristic is $2$. Beyond the commutative case, Question~\ref{q:main} has, to our knowledge, not been addressed, although $2\times2$ matrix rings over arbitrary rings arise naturally in several contexts. We describe four of them, which together explain the relevance of the question and the shape of its answer.

The first context is that of matrix units and inverse subsemigroups. A \emph{system of $2\times2$ matrix units} in a ring $R$ consists of nonzero elements $e_{11},e_{12},e_{21},e_{22}$ with $e_{ij}e_{kl}=\delta_{jk}e_{il}$; we do not require $e_{11}+e_{22}$ to be an identity. Such a system identifies a corner $eRe\cong\Mat_2(B)$, where $e=e_{11}+e_{22}$ and $B=e_{11}Re_{11}$ \cite[(17.5)]{Lam99}, and conversely every ring $\Mat_2(B)$ with $B$ unital contains one. The subring generated by the system is isomorphic to $\Mat_2(\Z/\ell\Z)$, where $\ell$ is the additive order of $e_{11}$, so that a system of matrix units splits into an arithmetic part, governed by $\Z/\ell\Z$, and a structural part, the corner ring $B$; both parts are visible in Question~\ref{q:main}. Already the arithmetic part allows non-Clifford phenomena in characteristic $2$: the ring $\Mat_2(\F_2)$ is Lie solvable of derived length $3$ and contains the rook monoid. Matrix units are also what produces non-Clifford inverse subsemigroups of the multiplicative semigroup $(R,\cdot)$. Indeed, $(R,\cdot)$ is an inverse semigroup only when it is a Clifford semigroup, whereas matrix rings contain non-Clifford inverse subsemigroups such as the rook monoid of partial permutation matrices \cite{So02}; the maximal inverse subsemigroups of the semigroup of linear transformations of a vector space were studied by Shneperman \cite{Sh74}. This leads to the second question that guides the paper.

\begin{question}\label{q:inverse}
How does Lie solvability of a ring $R$ restrict the inverse subsemigroups of its multiplicative semigroup $(R,\cdot)$?
\end{question}

We shall see that $(R,\cdot)$ has a non-Clifford inverse subsemigroup precisely when $R$ contains a system of $2\times2$ matrix units, and that systems of $n\times n$ matrix units with $n\ge3$ never occur in a Lie solvable ring; in particular, among the Brandt semigroups $B_n$ only $B_2$ can occur in its multiplicative semigroup. Question~\ref{q:inverse} is therefore a question about the corners $\Mat_2(B)$ of $R$, that is, a special case of Question~\ref{q:main}.

The second context is Jordan theory, and it explains why the answer involves Jordan products. The Lie ring $\gl_2(N)$, that is, $\Mat_2(N)$ under the commutator, carries the three-term grading
\[
\gl_2(N)=NE_{21}\oplus(NE_{11}\oplus NE_{22})\oplus NE_{12},
\]
and the brackets between its components of degree $\pm1$ produce the Jordan triple product of $N$:
\[
\bigl[[aE_{12},bE_{21}],cE_{12}\bigr]=(abc+cba)E_{12}\qquad(a,b,c\in N).
\]
This is the associative shadow of the Tits--Kantor--Koecher construction \cite[Section~0.3]{McC04}, in which the pair $(NE_{12},NE_{21})$ is a Jordan pair. When an identity is available, taking $c=1$ yields $(ab+ba)E_{12}$; thus the Lie structure of $\Mat_2(N)$ detects the Jordan product $x\circ y=xy+yx$ of $N$ rather than its commutator. The prime $2$ enters at precisely this point: for commutative $N$ one has $x\circ y=2xy$, whereas $x\circ y=[x,y]$ when $2N=0$. The same prime is visible in the grading itself. The elements $E_{12}$, $E_{21}$ and $h=E_{11}-E_{22}$ form an $\mathfrak{sl}_2$-triple, with $[h,E_{12}]=2E_{12}$ and $[h,E_{21}]=-2E_{21}$, and the grading is the decomposition of $\gl_2$ into the eigenspaces of $\operatorname{ad}h$ for the eigenvalues $2$, $0$ and $-2$. In characteristic $2$ these eigenvalues coincide, $h=E_{11}+E_{22}$ becomes central, and the triple degenerates into a Heisenberg algebra. This degeneration is the reason why $\gl_2(K)$ is solvable when $\chr K=2$. Over the integers the eigenvalues $\pm2$ survive as powers of $2$ accumulating along the derived series, and measuring this accumulation exactly (Section~\ref{s:prelim}) provides the arithmetic input for unital coefficient rings and matrix units.

The third context is structural. The ring $\Mat_2(B)$ is Morita equivalent to $B$, and most ring-theoretic properties are Morita invariant; Lie solvability is not, since $\Mat_2(B)$ is never Lie solvable for a nonzero commutative unital ring $B$ in which $2$ is invertible. Question~\ref{q:main} thus asks how Lie solvability behaves under the smallest nontrivial Morita equivalence $B\rightsquigarrow\Mat_2(B)$. It is also the case of Riley's problem on tensor products of Lie solvable algebras in which one factor is $\Mat_2(K)$: Riley \cite{Ri96} proved that $A\otimes B$ is Lie solvable when each factor is Lie nilpotent or Lie metabelian, whereas $\Mat_2(K)$ in characteristic $2$ has derived length $3$ and is not Lie nilpotent.

The fourth context concerns convolution algebras of groupoids. Leavitt path algebras and inverse semigroup algebras are Steinberg algebras of ample groupoids \cite{CFST14,CS15,St10}, and an orbit with two points contributes a block $\Mat_2(KH)$, where $H$ is the isotropy group; the simplest instance is $L_K(A_2)\cong\Mat_2(K)$ for the graph $A_2$ with one edge between two vertices. Similarly, the Leavitt path algebra of a cycle of length $2$ is $\Mat_2(K[x,x^{-1}])$, and a $\cD$-class with two idempotents and maximal subgroup $G$ contributes a block $\Mat_2(KG)$ to the algebra of a finite inverse semigroup \cite{St06}. Deciding when such algebras are Lie solvable therefore requires an understanding of $\Mat_2(KH)$ for arbitrary groups $H$, and the classical problem of Passi, Passman and Sehgal reappears here with an additional matrix factor.

\subsection*{Main results}
Our answer to Question~\ref{q:main} is formulated in terms of the \emph{Jordan ideal series} of $N$, defined by $\cJ_0(N)=N$, where $\cJ_{k+1}(N)$ is the ideal of $N$ generated by the Jordan products $x\circ y=xy+yx$ with $x,y\in\cJ_k(N)$, and of its length $\jl(N)$ (Definition~\ref{def:jordan}). Theorem~\ref{thm:main} shows that $\Mat_2(N)$ is Lie solvable if and only if this series reaches $0$, and that
\[
\jl(N)\le\dl(\Mat_2(N))\le3\jl(N),
\]
the constant $3$ being best possible. When $2N=0$ the condition becomes strong Lie solvability of $N$, and when $2$ is invertible on $N$ it becomes nilpotency of $N$; for $n\ge3$, the ring $\Mat_n(N)$ is Lie solvable exactly when $N$ is nilpotent (Proposition~\ref{prop:Mn}). For unital coefficient rings, the closed form of the derived series of $\gl_2(\Z)$ (Proposition~\ref{prop:gl2}) yields the exact derived lengths of the rings $\Mat_2(\Z/2^m\Z)$ and, with Theorem~\ref{thm:main}, a criterion for $\Mat_2(B)$ (Corollary~\ref{cor:M2B}). Question~\ref{q:inverse} is answered by Theorem~\ref{thm:units} and Corollary~\ref{cor:inverse}: in a Lie solvable ring of derived length $d$, the idempotent $e_{11}$ of a system of $2\times2$ matrix units has additive order $2^m$ with $m\le2\lfloor2^{d-1}/3\rfloor$, a sharp bound, and its corner ring $e_{11}Re_{11}$ is strongly Lie solvable. Finally, Corollary~\ref{cor:M2KG} and Theorem~\ref{thm:steinberg} characterize the Lie solvable algebras $\Mat_2(KG)$ in characteristic $2$ and the Lie solvable Steinberg algebras, with consequences for Leavitt path algebras and for algebras of finite inverse semigroups (Corollaries~\ref{cor:LPA} and~\ref{cor:finite-inverse}).

\subsection*{Organization of the paper}

Section~\ref{s:prelim} fixes notation, recalls the known results that we use, and proves the auxiliary results on which the main theorem rests. Section~\ref{s:main} proves the main theorem and treats matrices of size at least $3$. Section~\ref{s:unital} applies it to unital rings, matrix units, inverse subsemigroups and group algebras, and Section~\ref{s:steinberg} to Steinberg algebras, Leavitt path algebras and inverse semigroup algebras. Section~\ref{s:remarks} compares our results with the literature and discusses the sharpness of the bounds.

\section{Preliminaries and auxiliary results}\label{s:prelim}

This section fixes notation and collects what is used later: standard facts on derived series and matrix units, background on inverse semigroups, group algebras and ample groupoids, and the two auxiliary computations on which the main theorem rests, namely the derived series of $\gl_2(\Z)$ 
and the basic properties of the Jordan ideal series. Known results are quoted with references and without proofs. All rings are associative but not necessarily unital. For a ring $N$ we write $\hat N=\Z\oplus N$ for the ring obtained by adjoining an identity, $N^k$ for the additive subgroup generated by the products of $k$ elements of $N$, $xE_{ij}$ for the matrix with $x\in N$ in position $(i,j)$ and $0$ elsewhere, and $\gl_2(N)$ for the ring $\Mat_2(N)$ regarded as a Lie ring under the commutator $[x,y]=xy-yx$. The ring $N$ is \emph{nilpotent} if $N^t=0$ for some $t\ge1$, and the least such $t$ is its \emph{index of nilpotency} $t(N)$; it is \emph{nil} if each of its elements is nilpotent.

For a Lie ring $L$ let $L^{(0)}=L$, and let $L^{(k+1)}=[L^{(k)},L^{(k)}]$ be the additive subgroup generated by the elements $[a,b]$ with $a,b\in L^{(k)}$. A ring $R$ is \emph{Lie solvable of derived length $d$}, written $\dl(R)=d$, if $d$ is the least integer with $R^{(d)}=0$; if there is no such $d$, we write $\dl(R)=\infty$. The ring $R$ has \emph{no $2$-torsion} if $2x=0$ implies $x=0$, and $2$ is \emph{invertible on} $R$ if multiplication by $2$ is a bijection of $R$ onto itself. The \emph{strong Lie derived series} of $R$ is defined by $\Delta^{(0)}(R)=R$, where $\Delta^{(k+1)}(R)$ is the two-sided ideal of $R$ generated by $[\Delta^{(k)}(R),\Delta^{(k)}(R)]$; the ring $R$ is \emph{strongly Lie solvable} if $\Delta^{(k)}(R)=0$ for some $k$, and the least such $k$ is its \emph{strong Lie derived length} $\dlL(R)$. The ring $R$ is \emph{Lie nilpotent} if $[\cdots[[R,R],R],\dots,R]=0$ for some number of brackets, and \emph{Lie metabelian} if $\dl(R)\le2$. Clearly $R^{(k)}\subseteq\Delta^{(k)}(R)$, so that $\dl(R)\le\dlL(R)$, and a surjective ring homomorphism $\pi\colon R\to R'$ satisfies $\pi(\Delta^{(k)}(R))=\Delta^{(k)}(R')$, because it maps the ideal generated by a set onto the ideal generated by its image. The inclusion $R^{(k)}\subseteq\Delta^{(k)}(R)$ is far from an equality: for a field $K$ of characteristic $2$, the algebra $\Mat_2(K)$ is Lie solvable (Proposition~\ref{prop:gl2}(c) below) but not strongly Lie solvable, since $\Delta^{(1)}(\Mat_2(K))$ contains $[E_{11},E_{12}]=E_{12}$ and is therefore the whole algebra. Finally, for an integer $m\ge1$ we put
\[
d_m=\min\bigl\{k\ge2:\ 2\lfloor2^{k-1}/3\rfloor\ge m\bigr\};
\]
thus $d_1=d_2=3$, $d_3=d_4=4$, $d_m=5$ for $5\le m\le10$, and $d_m=6$ for $11\le m\le20$.

Matrix units are the bridge between arbitrary rings and the matrix rings studied in this paper. For $n\ge2$, a \emph{system of $n\times n$ matrix units} in a ring $R$ is a family $(e_{ij})_{1\le i,j\le n}$ of nonzero elements of $R$ such that $e_{ij}e_{kl}=\delta_{jk}e_{il}$ for all $i,j,k,l$; we do not require that $\sum_ie_{ii}=1$. By a classical theorem \cite[(17.5)]{Lam99}, applied to the ring $eRe$ with identity $e=\sum_ie_{ii}$, the element $e$ is an idempotent, $B=e_{11}Re_{11}$ is a ring with identity $e_{11}$, and
\begin{equation}\label{eq:corner}
eRe\longrightarrow\Mat_n(B),\qquad x\longmapsto(e_{1i}xe_{j1})_{i,j},
\end{equation}
is a ring isomorphism that maps $e_{ij}$ to $E_{ij}$, with inverse $(b_{ij})\mapsto\sum_{i,j}e_{i1}b_{ij}e_{1j}$. Since $ke_{ij}=e_{i1}(ke_{11})e_{1j}$ and $ke_{11}=e_{1i}(ke_{ij})e_{j1}$ for $k\in\Z$, all $e_{ij}$ have the same additive order $\ell\ge0$, where $\ell=0$ means infinite order, and $\ell\ne1$. Applied to the subring generated by the $e_{ij}$, which is the additive subgroup that they generate and has identity $e$, the isomorphism \eqref{eq:corner} shows that this subring is isomorphic to $\Mat_n(\Z/\ell\Z)$. A system of matrix units thus consists of an arithmetic skeleton $\Mat_n(\Z/\ell\Z)$ and a corner ring $B$.

We shall constantly use the following elementary properties of derived series, in which $L$ is a Lie ring and $R$ is a ring.
\begin{enumerate}[label=\textup{(D\arabic*)},ref=(D\arabic*),leftmargin=3em]
\item\label{D:sub} If $L'$ is a Lie subring of $L$, then $L'^{(k)}\subseteq L^{(k)}$ for all $k$.
\item\label{D:hom} If $\psi\colon L\to L''$ is a surjective homomorphism of Lie rings, then $\psi(L^{(k)})=L''^{(k)}$ for all $k$.
\item\label{D:subdirect} If $R$ is a \emph{subdirect product} of rings $R_i$, that is, $R\subseteq\prod_iR_i$ and every projection $R\to R_i$ is surjective, then $\dl(R)=\sup_i\dl(R_i)$.
\item\label{D:ideal} If $J$ is an ideal of $R$, then $J^{(k)}\subseteq J^{2^k}$ for all $k$; in particular, $\dl(J)\le\lceil\log_2t\rceil$ if $J^t=0$.
\item\label{D:units} If $R$ contains a system of $n\times n$ matrix units with $n\ge3$, then $R$ is not Lie solvable.
\end{enumerate}
Properties \ref{D:sub}--\ref{D:subdirect} follow directly from the definitions, and \ref{D:ideal} follows by induction on $k$ from the inclusion $[J^a,J^b]\subseteq J^{a+b}$. Property \ref{D:units} is classical, and we recall its proof: for pairwise distinct $i,j,k$ we have $[e_{ij},e_{jk}]=e_{ik}$, so if all $e_{ik}$ with $i\ne k$ lie in $R^{(r)}$, then, choosing $j\notin\{i,k\}$, they also lie in $R^{(r+1)}$; since $e_{12}\ne0$, no term of the derived series vanishes.

Matrix units also govern the inverse subsemigroups of the multiplicative semigroup $(R,\cdot)$. A semigroup $S$ is \emph{inverse} if every $s\in S$ has a unique $s^{-1}\in S$ with $ss^{-1}s=s$ and $s^{-1}ss^{-1}=s^{-1}$; equivalently, $S$ is regular and its idempotents commute \cite[Theorem~5.1.1]{How95}. Put $\dd(s)=s^{-1}s$ and $\rr(s)=ss^{-1}$. An inverse semigroup is a \emph{Clifford semigroup} if its idempotents are central, equivalently if $\dd(s)=\rr(s)$ for all $s$ \cite[Theorem~4.2.1]{How95}. Two elements $s,t$ of an inverse semigroup $S$ are \emph{$\cD$-related} if $\dd(s)=\dd(u)$ and $\rr(u)=\rr(t)$ for some $u\in S$; the equivalence classes of this relation are the \emph{$\cD$-classes} of $S$, and the \emph{maximal subgroup} of $S$ at an idempotent $e$ is $\{s\in S:\dd(s)=\rr(s)=e\}$ \cite[Corollary~2.2.6 and Section~5.1]{How95}. For $n\ge2$, the \emph{Brandt semigroup} $B_n$ consists of a zero $0$ and elements $e_{ij}$ with $1\le i,j\le n$, with product $e_{ij}e_{kl}=e_{il}$ if $j=k$ and $e_{ij}e_{kl}=0$ otherwise; more generally, for a group $G$, the Brandt semigroup $B(G,n)$ consists of $0$ and the triples $(i,g,j)$ with $1\le i,j\le n$ and $g\in G$, with product $(i,g,j)(k,h,l)=(i,gh,l)$ if $j=k$ and $0$ otherwise. The semigroup $B_2$ and the rook monoid $R_2$ of $2\times2$ partial permutation matrices \cite[Section~1]{So02} are inverse but not Clifford, since $\dd(e_{12})=e_{22}\ne e_{11}=\rr(e_{12})$. By contrast, if $(R,\cdot)$ is itself inverse, then its idempotents commute and are therefore central: for an idempotent $e$ and $r\in R$, the elements $e+er-ere$ and $e+re-ere$ are idempotents, and since they commute with $e$ we obtain $er=ere=re$. Hence $(R,\cdot)$ is then a Clifford semigroup.

The classical results on which our answer to Question~\ref{q:main} builds are collected in the following theorem; part (a) is the commutative case of Question~\ref{q:main}, and part (b) is needed for the applications to group algebras. Recall that, for a prime $p$, a group $G$ is \emph{$p$-abelian} if its commutator subgroup $G'$ is a finite $p$-group, and that $0$-abelian means abelian.

\begin{theorem}\label{thm:PPS}
Let $K$ be a field of characteristic $p\ge0$.
\begin{enumerate}[label=\textup{(\alph*)}]
\item \cite[Lemma~1.2]{PPS73} If $S$ is a commutative $K$-algebra, not necessarily unital, which is not nilpotent, then $\Mat_n(S)$ is Lie solvable if and only if $n=1$, or $n=2$ and $p=2$.
\item \cite[Theorem, p.~748]{PPS73} The group algebra $KG$ is Lie solvable if and only if either $p\ne2$ and $G$ is $p$-abelian, or $p=2$ and $G$ has a $2$-abelian subgroup of index at most $2$.
\end{enumerate}
\end{theorem}

Let $K$ be a field of characteristic $p>0$ and $G$ a group, and denote by $\omega(KG)$ the augmentation ideal of $KG$, that is, the kernel of the augmentation map $KG\to K$, $\sum_ga_gg\mapsto\sum_ga_g$. If $H$ is a normal subgroup of $G$, then $\omega(KH)KG=KG\,\omega(KH)$ is the kernel of the natural map $KG\to K[G/H]$, and $(\omega(KH)KG)^t=\omega(KH)^tKG$ for all $t\ge1$. If $\omega(KH)$ is nilpotent, we write $t(H)=t(\omega(KH))$ for its index of nilpotency; thus $t(1)=1$. Explicit bounds for group algebras will be derived from the following known facts.

\begin{theorem}\label{thm:groupalg}
Let $K$ be a field of characteristic $p>0$ and $G$ a group.
\begin{enumerate}[label=\textup{(\alph*)}]
\item \cite[Theorem~3.7]{Je41} If $H$ is a finite $p$-group, then $\omega(KH)$ is nilpotent and $1+(p-1)\log_p|H|\le t(H)\le|H|$.
\item \cite[Theorem~V.5.1]{Seh78} The algebra $KG$ is strongly Lie solvable if and only if $G'$ is a finite $p$-group.
\item \cite{CS07} (see also \cite[p.~295]{Sp08}) If $KG$ is strongly Lie solvable and not commutative, then $\dlL(KG)\ge\lceil\log_2(t(G')+1)\rceil$.
\end{enumerate}
\end{theorem}

The lower bound in (a) follows from the formula $t(H)=1+(p-1)\sum_{i\ge1}i\,d_i$ of \cite[Theorem~3.7]{Je41}, where $p^{d_i}=|D_i:D_{i+1}|$ for the dimension subgroups $D_i=\{h\in H:h-1\in\omega(KH)^i\}$ of $H$, since $\sum_id_i=\log_p|H|$; the upper bound follows from $\dim_K\omega(KH)=|H|-1$, because the powers of a nilpotent ideal decrease strictly until they vanish.

Leavitt path algebras and inverse semigroup algebras are convolution algebras of ample groupoids, and we recall the relevant notions. Following \cite[Definition~2.4]{St10}, a topological groupoid $\cG$ is \emph{ample} if its unit space $\cG^{(0)}$ is locally compact and Hausdorff, its source map $s$, and hence its range map $r$, is a local homeomorphism, and its compact open \emph{bisections}, the open subsets $U$ on which $s$ and $r$ are injective, form a basis of its topology. We do not assume that $\cG$ is Hausdorff; for Steinberg algebras of non-Hausdorff groupoids see also \cite[Section~3]{CEPSS19}. The compact open bisections form an inverse semigroup under $UV=\{\alpha\beta:\alpha\in U,\ \beta\in V,\ s(\alpha)=r(\beta)\}$ and $U^{-1}=\{\alpha^{-1}:\alpha\in U\}$, and every compact open subset of $\cG^{(0)}$ is one of them. The \emph{Steinberg algebra} $A_K(\cG)$ over a field $K$ is the subspace of $K^{\cG}$ spanned by the characteristic functions $1_U$ of compact open bisections \cite[Proposition~3.3]{St10}, with the convolution $(f*g)(\gamma)=\sum_{\alpha\beta=\gamma}f(\alpha)g(\beta)$; one has $1_U*1_V=1_{UV}$ \cite[Proposition~3.5]{St10}. The \emph{orbit} of $x\in\cG^{(0)}$ is $\{r(\gamma):s(\gamma)=x\}$, the \emph{isotropy group} at $x$ is $\cG_x^x=\{\gamma:s(\gamma)=r(\gamma)=x\}$, and $\cG$ is a \emph{group bundle} if every orbit is a single point. For an orbit $O$, the set $\cG|_O=s^{-1}(O)$ is a groupoid with unit space $O$ and a single orbit, that is, a \emph{transitive} groupoid, and $\cG$ is the disjoint union of the groupoids $\cG|_O$. For a discrete groupoid $\mathcal H$, $K[\mathcal H]$ denotes the algebra of finitely supported functions on $\mathcal H$ with the same convolution, and $\delta_\gamma$ the characteristic function of $\{\gamma\}$. If $\mathcal H$ is transitive with finite unit space $\{x_1,\dots,x_n\}$ and isotropy group $H$ at $x_1$, a choice of arrows $\gamma_i$ from $x_1$ to $x_i$ writes every arrow from $x_j$ to $x_i$ uniquely as $\gamma_ih\gamma_j^{-1}$ with $h\in H$, and $\delta_{\gamma_ih\gamma_j^{-1}}\mapsto hE_{ij}$ is an algebra isomorphism $K[\mathcal H]\cong\Mat_n(KH)$.

When all orbits are finite, the Steinberg algebra is controlled by matrix algebras over the isotropy group algebras. For a finite unit space this is the decomposition $A_K(\cG)\cong\bigoplus_O\Mat_{|O|}(KH_O)$ of \cite[Proposition~3.1]{St18} (see also \cite[Proposition~4.3]{CH20}). In Section~\ref{s:steinberg} we need the following extension to infinite unit spaces, where the direct sum is replaced by a subdirect product. Its proof restricts functions to the finite orbits, which are closed invariant subsets of the unit space; we include a proof, since we could not find this statement in the literature.

\begin{proposition}\label{prop:finite-orbits}
Let $\cG$ be an ample groupoid all of whose orbits are finite, and let $K$ be a field. For each orbit $O$ choose $x_O\in O$ and put $H_O=\cG_{x_O}^{x_O}$. Then restriction $f\mapsto f|_{\cG|_O}$ is a surjective algebra homomorphism $A_K(\cG)\to K[\cG|_O]\cong\Mat_{|O|}(KH_O)$ for each $O$, and these homomorphisms together embed $A_K(\cG)$ as a subdirect product of the algebras $\Mat_{|O|}(KH_O)$. In particular, $\dl(A_K(\cG))=\sup_O\dl\bigl(\Mat_{|O|}(KH_O)\bigr)$.
\end{proposition}

\begin{proof}
Let $O$ be an orbit. A compact open bisection meets each fiber $s^{-1}(x)$ in at most one point, and $\cG|_O=\bigcup_{x\in O}s^{-1}(x)$ with $O$ finite; hence the restriction of every $f\in A_K(\cG)$ to $\cG|_O$ has finite support. If $\alpha\beta=\gamma\in\cG|_O$, then $s(\beta)=s(\gamma)\in O$ and $s(\alpha)=r(\beta)\in O$, so $\alpha,\beta\in\cG|_O$; consequently $(f*g)|_{\cG|_O}=f|_{\cG|_O}*g|_{\cG|_O}$, and restriction is a homomorphism $\rho_O$. It is surjective: for $\gamma\in\cG|_O$, choose a compact open bisection $U\ni\gamma$ and a compact open $V\subseteq\cG^{(0)}$ with $V\cap O=\{s(\gamma)\}$; then $UV$ is a compact open bisection with $UV\cap\cG|_O=\{\gamma\}$, so $\rho_O(1_{UV})=\delta_\gamma$. Since $\cG$ is the disjoint union of the $\cG|_O$, the family $(\rho_O)_O$ is injective. The last assertion is \ref{D:subdirect}.
\end{proof}

We now turn to the auxiliary results that enter the proof of the main theorem and its applications. The arithmetic of matrix units is governed by the rings $\Mat_2(\Z/\ell\Z)$, and matrix rings over commutative rings annihilated by a power of $2$ enter all the bounds for unital coefficient rings. Both are controlled by the derived series of the Lie ring $\gl_2(\Z)$, which the following proposition determines in closed form; the factors $2$ in the relations $[h,x]=2x$ and $[h,y]=-2y$ are responsible for the exponents that appear. For a commutative algebra $S$ of characteristic $2$, the underlying computation appears in the proof of \cite[Lemma~1.2]{PPS73}, which shows that the second derived term of $\Mat_2(S)$ consists of scalar matrices; the proposition refines it by keeping track of the powers of $2$.

\begin{proposition}\label{prop:gl2}
Let $L$ denote the Lie ring $\Mat_2(\Z)$, and put $x=E_{12}$, $y=E_{21}$ and $h=E_{11}-E_{22}$.
\begin{enumerate}[label=\textup{(\alph*)}]
\item For every $k\ge1$,
\[
L^{(k)}=\Z\,2^{\alpha_k}h+\Z\,2^{\beta_k}x+\Z\,2^{\beta_k}y,\qquad\text{where}\quad \alpha_k=2\lfloor2^{k-1}/3\rfloor\quad\text{and}\quad\beta_k=\lfloor2^k/3\rfloor .
\]
Moreover, $\alpha_k\le\beta_k\le\alpha_k+1$ for all $k$, and the sequence $(\alpha_k)_{k\ge2}$ is strictly increasing.
\item Let $\ell\ge0$ be an integer with $\ell\ne1$. Then $\Mat_2(\Z/\ell\Z)$ is Lie solvable if and only if $\ell=2^m$ with $m\ge1$, and in this case $\dl(\Mat_2(\Z/2^m\Z))=d_m=\min\{k\ge2:\alpha_k\ge m\}$.
\item If $A$ is a commutative ring, not necessarily unital, with $2^mA=0$ for some $m\ge1$, then $\dl(\Mat_2(A))\le d_m$, with equality if $A$ is unital of characteristic exactly $2^m$. In particular, $\dl(\Mat_2(A))=3$ for every unital commutative ring $A$ of characteristic $2$.
\end{enumerate}
\end{proposition}

\begin{proof}
(a) We have $[x,y]=h$, $[h,x]=2x$ and $[h,y]=-2y$. The group $L^{(1)}$ is generated by the commutators of the matrix units: $[E_{11},E_{12}]=x$, $[E_{11},E_{21}]=-y$, $[E_{12},E_{21}]=h$, $[E_{12},E_{22}]=x$, $[E_{21},E_{22}]=-y$ and $[E_{11},E_{22}]=0$. Hence $L^{(1)}=\Z h+\Z x+\Z y$, which is the assertion for $k=1$, as $\alpha_1=\beta_1=0$. Define integers $a_k,b_k$ by $a_1=b_1=0$, $a_{k+1}=2b_k$ and $b_{k+1}=a_k+b_k+1$. If $L^{(k)}=\Z\,2^{a_k}h+\Z\,2^{b_k}x+\Z\,2^{b_k}y$, then, since the Lie product is bilinear, $L^{(k+1)}$ is generated by the commutators
\[
[2^{b_k}x,2^{b_k}y]=2^{2b_k}h,\qquad
[2^{a_k}h,2^{b_k}x]=2^{a_k+b_k+1}x,\qquad
[2^{a_k}h,2^{b_k}y]=-2^{a_k+b_k+1}y,
\]
so the description holds for all $k$ with exponents $a_k,b_k$. By induction, $a_k\le b_k\le a_k+1$: this holds for $k=1$, and it implies $a_{k+1}=2b_k\le a_k+b_k+1=b_{k+1}\le2b_k+1=a_{k+1}+1$. It remains to identify the exponents. We have $b_1=0$, $b_2=1$ and $b_{k+1}=2b_{k-1}+b_k+1$ for $k\ge2$. Put $\beta'_k=(2^k-1)/3$ for even $k$ and $\beta'_k=(2^k-2)/3$ for odd $k$, so that $\beta'_k=\lfloor2^k/3\rfloor$. For even $k$ we obtain $2\beta'_{k-1}+\beta'_k+1=\bigl(2(2^{k-1}-2)+2^k-1+3\bigr)/3=(2^{k+1}-2)/3=\beta'_{k+1}$, and for odd $k$ we obtain $2\beta'_{k-1}+\beta'_k+1=\bigl(2(2^{k-1}-1)+2^k-2+3\bigr)/3=(2^{k+1}-1)/3=\beta'_{k+1}$. Since $\beta'_1=0$ and $\beta'_2=1$, it follows that $b_k=\beta_k$ for all $k$, hence $a_k=2\beta_{k-1}=\alpha_k$ for $k\ge2$, and $a_1=0=\alpha_1$. Finally, $\alpha_2=0<2=\alpha_3$, and $\alpha_{k+1}=2\beta_k\ge2\alpha_k>\alpha_k$ for $k\ge3$.

(b) The projection $\Z\to\Z/\ell\Z$ induces a surjective homomorphism of Lie rings $\Mat_2(\Z)\to\Mat_2(\Z/\ell\Z)$, so by \ref{D:hom} the derived series of $\Mat_2(\Z/\ell\Z)$ is the image of the series in (a). If $\ell=0$, then $2^{\beta_k}x\ne0$ lies in $L^{(k)}$ for all $k$. If $p$ is an odd prime divisor of $\ell$, then the image of $L^{(k)}$ in $\Mat_2(\F_p)$ contains $2^{\beta_k}x\ne0$, since $2$ is invertible in $\F_p$, and \ref{D:hom} shows that $\Mat_2(\Z/\ell\Z)$ is not Lie solvable. If $\ell=2^m$, then $h,x,y$ are linearly independent over $\Z/2^m\Z$, so the image of $L^{(k)}$ vanishes if and only if $2^m$ divides $2^{\alpha_k}$ and $2^{\beta_k}$, that is, $\alpha_k\ge m$; since $\alpha_1=\alpha_2=0<m$, the derived length is $\min\{k\ge2:\alpha_k\ge m\}=d_m$.

(c) For $X,Y\in\Mat_2(\Z)$ and $a,b\in A$ we have $[aX,bY]=ab[X,Y]$, because $A$ is commutative. By induction on $k$, the group $\Mat_2(A)^{(k)}$ is contained in the additive subgroup generated by the elements $aZ$ with $a\in A$ and $Z\in L^{(k)}$, that is, in $A\,2^{\alpha_k}h+A\,2^{\beta_k}x+A\,2^{\beta_k}y$ by (a). For $k=d_m$ we have $\beta_k\ge\alpha_k\ge m$, so this group is $0$. If $A$ is unital of characteristic $2^m$, then $\Mat_2(\Z/2^m\Z)\cong\Mat_2(\Z\cdot1_A)$ is a subring of $\Mat_2(A)$, and \ref{D:sub} and (b) give $\dl(\Mat_2(A))\ge d_m$.
\end{proof}

Table~\ref{tab:exp} lists the first exponents. The closed form gives $\alpha_k=(2^k-2)/3$ for odd $k$ and $\alpha_k=(2^k-4)/3$ for even $k\ge2$. Since $(2^k-4)/3\le\alpha_k<2^k/3$, we obtain $3m<2^{d_m}$ and $2^{d_m-1}<3m+4$, so that $d_m$ equals $\lfloor\log_2(3m)\rfloor+1$ or $\lfloor\log_2(3m)\rfloor+2$. For a field $K$ of characteristic $2$, explicitly, $h=E_{11}+E_{22}=1$ and the derived series of $\gl_2(K)$ is $\gl_2(K)\supset\spn\{x,y,1\}\supset K\cdot1\supset0$.

\begin{table}[ht]
\centering
\begin{tabular}{@{}c|cccccccc@{}}
\toprule
$k$ & $1$ & $2$ & $3$ & $4$ & $5$ & $6$ & $7$ & $8$\\
\midrule
$\alpha_k$ & $0$ & $0$ & $2$ & $4$ & $10$ & $20$ & $42$ & $84$\\
$\beta_k$ & $0$ & $1$ & $2$ & $5$ & $10$ & $21$ & $42$ & $85$\\
\bottomrule
\end{tabular}
\caption{The exponents in $L^{(k)}=\Z\,2^{\alpha_k}h+\Z\,2^{\beta_k}x+\Z\,2^{\beta_k}y$.}
\label{tab:exp}
\end{table}

The second ingredient of the main theorem is the three-term grading of $\gl_2(N)$, which links the Lie structure of $\Mat_2(N)$ to the Jordan structure of $N$. Let $N$ be a ring, $M=\Mat_2(N)$ and $\hat M=\Mat_2(\hat N)$; thus $M$ is an ideal of the unital ring $\hat M$. Put $M_1=NE_{12}$, $M_0=NE_{11}\oplus NE_{22}$ and $M_{-1}=NE_{21}$. Direct multiplication gives, for $a,a',b,c\in\hat N$,
\begin{equation}\label{eq:grading}
\begin{gathered}
[aE_{12},bE_{21}]=abE_{11}-baE_{22},\\
[aE_{11}+a'E_{22},\,cE_{12}]=(ac-ca')E_{12},\qquad [aE_{11}+a'E_{22},\,cE_{21}]=(a'c-ca)E_{21},
\end{gathered}
\end{equation}
and, combining the first two identities,
\begin{equation}\label{eq:triple}
\bigl[[aE_{12},bE_{21}],cE_{12}\bigr]=(abc+cba)E_{12},
\end{equation}
which for $c=1$ shows that $\bigl[[aE_{12},bE_{21}],E_{12}\bigr]=(ab+ba)E_{12}$ for all $a,b\in\hat N$. Together with $[M_1,M_1]=[M_{-1},M_{-1}]=0$, the identities \eqref{eq:grading} show that $M=M_{-1}\oplus M_0\oplus M_1$ is a grading of the Lie ring $M$. Finally, $P=E_{12}+E_{21}$ satisfies $P^2=1$ and $P(aE_{12})P=aE_{21}$, so that $Y\mapsto PYP$ is an automorphism of $\hat M$ that maps $M$ onto $M$ and interchanges $M_1$ and $M_{-1}$. The next remark places these identities in their natural context and explains the form of the invariant introduced below.

\begin{remark}\label{rem:TKK}
The identities \eqref{eq:grading} and \eqref{eq:triple} are the associative shadow of the Tits--Kantor--Koecher construction \cite[Section~0.3]{McC04}: the pair $(M_1,M_{-1})\cong(N,N)$, with the triple products $\{a,b,c\}=abc+cba$ of \eqref{eq:triple}, is a Jordan pair \cite[Definition~1.2]{Loo75}, and $M_0$ acts on it through \eqref{eq:grading}. Thus the Lie structure of $\Mat_2(N)$ detects the Jordan triple product of $N$ and, once an identity is available as in $\hat M$, the Jordan product $a\circ b=ab+ba$. This is why the invariant that controls the Lie solvability of $\Mat_2(N)$ is defined by Jordan products and not by commutators.
\end{remark}

We now introduce that invariant; the terminology should not be confused with the Jordan ideals of \cite[Chapter~1]{Her69}, which are the additive subgroups $U$ of a ring $R$ with $u\circ x\in U$ for all $u\in U$ and $x\in R$. Since the main theorem concerns rings without identity, the series is formed with two-sided ideals.

\begin{definition}\label{def:jordan}
Let $N$ be a ring and put $x\circ y=xy+yx$. The \emph{Jordan ideal series} of $N$ is defined by $\cJ_0(N)=N$, where $\cJ_{k+1}(N)$ is the two-sided ideal of $N$ generated by the elements $x\circ y$ with $x,y\in\cJ_k(N)$; it is a descending series of ideals. We put $\jl(N)=\min\{k:\cJ_k(N)=0\}$, and $\jl(N)=\infty$ if there is no such $k$. A ring $A$ is \emph{anticommutative} if $x\circ y=0$ for all $x,y\in A$.
\end{definition}

The Jordan ideal series is the ideal-theoretic counterpart of the derived series of the Jordan ring of $N$, just as the strong Lie derived series is that of the derived series of its Lie ring. The next remark compares the two notions and explains why the ideal-theoretic version is the relevant one.

\begin{remark}\label{rem:Nplus}
Let $N^+$ be the additive group of $N$ with the product $x\circ y$, and define the derived series of this Jordan ring by $(N^+)^{(0)}=N$ and $(N^+)^{(k+1)}=(N^+)^{(k)}\circ(N^+)^{(k)}$, the additive subgroup generated by the products $x\circ y$ with $x,y\in(N^+)^{(k)}$. Induction gives $(N^+)^{(k)}\subseteq\cJ_k(N)$, so $N^+$ is solvable of length at most $\jl(N)$, and hence at most $\dl(\Mat_2(N))$ by Theorem~\ref{thm:main} below. The converse fails: for a field $K$ of characteristic $2$, the Jordan ring of $KS_3$ coincides with its Lie ring, which is solvable by Theorem~\ref{thm:PPS}(b), whereas $\Mat_2(KS_3)$ is not Lie solvable (Example~\ref{ex:tensor}). The ideal-theoretic series appears because every term of the derived series of $\Mat_2(N)$ is stable under the brackets with $NE_{11}$ and $NE_{22}$, which multiply the entries of $M_1$ by elements of $N$ (see the proof of Theorem~\ref{thm:main}).
\end{remark}

By definition, the quotients $\cJ_k(N)/\cJ_{k+1}(N)$ are anticommutative, and the upper bound in the main theorem rests on the following observation: in an anticommutative ring all products are central and all triple products are annihilated by $2$, so that the second derived term of the $2\times2$ matrix ring consists of central scalar matrices. Since a commutative ring annihilated by $2$ is anticommutative, the lemma extends the computation in the proof of \cite[Lemma~1.2]{PPS73}.

\begin{lemma}\label{lem:anticomm}
Let $A$ be an anticommutative ring. Then $A^2$ is contained in the center of $A$, $2A^3=0$, and
\[
\Mat_2(A)^{(1)}\subseteq A^2E_{12}+A^2E_{21}+A^2\cdot1+2A^2E_{11}+2A^2E_{22},\qquad \Mat_2(A)^{(2)}\subseteq\{d\cdot1:d\in A^4\}.
\]
In particular, $\Mat_2(A)^{(2)}$ is central in $\Mat_2(A)$, and $\dl(\Mat_2(A))\le3$. If $A$ has no $2$-torsion, then $A^3=0$ and $\dl(\Mat_2(A))\le2$.
\end{lemma}

\begin{proof}
For $x,y,z\in A$, repeated use of $uv=-vu$ gives $(xy)z=-z(xy)=-(zx)y=(xz)y=x(zy)=-x(yz)$, so that $2xyz=0$ and $[xy,z]=(xy)z-z(xy)=2xyz=0$. Hence $A^2$ is central in $A$ and $2A^3=0$. For $a,b\in A$ and integer matrices $X,Y$ we have $[aX,bY]=abXY-baYX=ab(X\circ Y)$, because $ba=-ab$, and the Jordan products of the matrix units are $0$, $2E_{11}$, $2E_{22}$, $E_{12}$, $E_{21}$ or $1$; for instance, $E_{11}\circ E_{12}=E_{12}$ and $E_{12}\circ E_{21}=1$. This gives the first inclusion. For $c,c'\in A^2$ we have $cc'=c'c$, since $c$ is central, and $2cc'=0$, since $cc'\in A^4\subseteq A^3$; hence $[cX,c'Y]=cc'[X,Y]$ for integer matrices $X,Y$. The commutators of the matrices $E_{12}$, $E_{21}$, $1$, $2E_{11}$ and $2E_{22}$ are integral combinations of $h=E_{11}-E_{22}$, $2E_{12}$ and $2E_{21}$, so $\Mat_2(A)^{(2)}$ is contained in $\{dh:d\in A^4\}$, which equals $\{d\cdot1:d\in A^4\}$ because $2A^4=0$. These scalar matrices are central in $\Mat_2(A)$, since $A^4\subseteq A^2$ is central in $A$, and therefore $\Mat_2(A)^{(3)}=0$. If $A$ has no $2$-torsion, then $A^3=0$, so $A^4=0$ and $\Mat_2(A)^{(2)}=0$.
\end{proof}

The last auxiliary result collects the properties of the Jordan ideal series needed later; parts (b) and (c) show how the series reflects the arithmetic of $2$ in $N$. Parts (a) and (b) follow directly from the definitions, and part (c) is obtained by applying Lemma~\ref{lem:anticomm} to the factors of the series.

\begin{lemma}\label{lem:jordan-series}
Let $N$ be a ring.
\begin{enumerate}[label=\textup{(\alph*)}]
\item A surjective ring homomorphism $\pi\colon N\to N'$ satisfies $\pi(\cJ_k(N))=\cJ_k(N')$ for all $k$. Moreover, $\cJ_{k+1}(N)\subseteq\cJ_k(N)^2$, so that $\cJ_k(N)\subseteq N^{2^k}$ for all $k$; in particular, $\jl(N)\le\lceil\log_2t(N)\rceil$ if $N$ is nilpotent.
\item If $2N=0$, then $\cJ_k(N)=\Delta^{(k)}(N)$ for all $k$, and $\jl(N)=\dlL(N)$.
\item Suppose that $2$ is invertible on $N$. Then every quotient $\cJ_k(N)/\cJ_{k+1}(N)$ has no $2$-torsion, $\cJ_k(N)^3\subseteq\cJ_{k+1}(N)$ and $N^{3^k}\subseteq\cJ_k(N)$ for all $k$. Hence $\jl(N)<\infty$ if and only if $N$ is nilpotent, and then $\log_3t(N)\le\jl(N)\le\lceil\log_2t(N)\rceil$. If, moreover, $N$ is commutative, then $\cJ_k(N)=N^{2^k}$ for all $k$.
\end{enumerate}
\end{lemma}

\begin{proof}
(a) The first assertion holds because $\pi$ maps the ideal generated by a set onto the ideal generated by its image. Since $x\circ y\in\cJ_k(N)^2$ for $x,y\in\cJ_k(N)$, and $\cJ_k(N)^2$ is an ideal, we get $\cJ_{k+1}(N)\subseteq\cJ_k(N)^2$, and induction gives $\cJ_k(N)\subseteq N^{2^k}$.

(b) If $2N=0$, then $x\circ y=[x,y]$, and the two series are defined by the same recursion.

(c) Write $\cJ_k=\cJ_k(N)$. By induction on $k$, each $\cJ_k$ is closed under division by $2$: if this holds for $\cJ_k$, then the set of products $x\circ y$ with $x,y\in\cJ_k$ is closed under division by $2$, because $(\frac12x)\circ y=\frac12(x\circ y)$, and hence so is the ideal $\cJ_{k+1}$ that it generates. Consequently, if $x\in\cJ_k$ and $2x\in\cJ_{k+1}$, then $x=\frac12(2x)\in\cJ_{k+1}$, so $\cJ_k/\cJ_{k+1}$ has no $2$-torsion. Lemma~\ref{lem:anticomm} now gives $\cJ_k^3\subseteq\cJ_{k+1}$, and induction gives $N^{3^k}\subseteq\cJ_k$, since $N^{3^{k+1}}=(N^{3^k})^3\subseteq\cJ_k^3\subseteq\cJ_{k+1}$. Thus $\jl(N)<\infty$ implies $t(N)\le3^{\jl(N)}$, and the remaining inequality is (a). If $N$ is commutative, then $x\circ y=2xy$, so $\cJ_{k+1}$ is the ideal generated by $2\cJ_k^2$, which equals $\cJ_k^2$; hence $\cJ_k=N^{2^k}$.
\end{proof}

\section{Lie solvability of \texorpdfstring{$\Mat_2(N)$}{M2(N)}}\label{s:main}

We now answer Question~\ref{q:main}. The lower bound is obtained by following the elements of degree $1$ through the derived series of $\Mat_2(N)$ by means of the identity \eqref{eq:triple} with $c=1$, and the upper bound by applying Lemma~\ref{lem:anticomm} to the anticommutative quotients of the Jordan ideal series. Part (3) contains, in particular, the case of algebras over fields of characteristic different from $2$, and at the end of the section we compare the answer with the case of matrices of size at least $3$.

\begin{theorem}\label{thm:main}
Let $N$ be a ring, not necessarily unital.
\begin{enumerate}[label=\textup{(\arabic*)}]
\item $\jl(N)\le\dl(\Mat_2(N))\le3\jl(N)$, and the constant $3$ is best possible. In particular, $\Mat_2(N)$ is Lie solvable if and only if $\cJ_k(N)=0$ for some $k$.
\item If $2N=0$, then $\dlL(N)\le\dl(\Mat_2(N))\le3\dlL(N)$. In particular, $\Mat_2(N)$ is Lie solvable if and only if $N$ is strongly Lie solvable.
\item If $2$ is invertible on $N$, for instance if $N$ is an algebra over a field of characteristic different from $2$, then $\Mat_2(N)$ is Lie solvable if and only if $N$ is nilpotent, and in this case
\[
\log_3t(N)\le\jl(N)\le\dl(\Mat_2(N))\le\lceil\log_2t(N)\rceil\le\lceil\jl(N)\log_23\rceil;
\]
if, moreover, $N$ is commutative, then $\dl(\Mat_2(N))=\lceil\log_2t(N)\rceil$.
\item If $N\ne0$ is unital, then $\jl(N)+1\le\dl(\Mat_2(N))$.
\end{enumerate}
\end{theorem}

\begin{proof}
We first prove the lower bounds in (1) and (4). Let $M=\Mat_2(N)$, $\hat M=\Mat_2(\hat N)$ and $P=E_{12}+E_{21}$ be as in Section~\ref{s:prelim}, and put $X_k=\{x\in N:xE_{12}\in M^{(k)}\}$. Since $M$ is an ideal of $\hat M$, we have $[M,\hat M]\subseteq M$; since $[A,B]$ is a Lie ideal of $\hat M$ whenever $A$ and $B$ are, by the Jacobi identity, every $M^{(k)}$ is a Lie ideal of $\hat M$. Moreover, every $M^{(k)}$ is stable under the automorphism $Y\mapsto PYP$, which preserves $M$. For $b\in N$ and $x\in X_k$, the identities \eqref{eq:grading} give $bxE_{12}=[bE_{11},xE_{12}]$ and $xbE_{12}=-[bE_{22},xE_{12}]$, so $X_k$ is an ideal of $N$. For $x,y\in X_k$ we have $yE_{21}=P(yE_{12})P\in M^{(k)}$, hence $[xE_{12},yE_{21}]\in M^{(k+1)}$, and \eqref{eq:triple} with $c=1$ gives
\[
(x\circ y)E_{12}=\bigl[[xE_{12},yE_{21}],\,E_{12}\bigr]\in[M^{(k+1)},\hat M]\subseteq M^{(k+1)}.
\]
Thus $x\circ y\in X_{k+1}$ for all $x,y\in X_k$. Since $X_{k+1}$ is an ideal and $X_0=N=\cJ_0(N)$, induction gives $\cJ_k(N)\subseteq X_k$ for all $k$. If $N$ is unital, then $X_1=N$, because $bE_{12}=[E_{11},bE_{12}]$, and the same induction gives $\cJ_k(N)\subseteq X_{k+1}$. As $M^{(d)}=0$ forces $X_d=0$, we obtain $\jl(N)\le\dl(M)$, and $\jl(N)+1\le\dl(M)$ if $N\ne0$ is unital.

For the upper bound in (1), write $\cJ_k=\cJ_k(N)$. The quotient $A_k=\cJ_k/\cJ_{k+1}$ is anticommutative, so $\dl(\Mat_2(A_k))\le3$ by Lemma~\ref{lem:anticomm}. Applying \ref{D:hom} to the projection $\Mat_2(\cJ_k)\to\Mat_2(A_k)$, we obtain $\Mat_2(\cJ_k)^{(3)}\subseteq\Mat_2(\cJ_{k+1})$. By induction and \ref{D:sub}, $\Mat_2(N)^{(3k)}\subseteq\Mat_2(\cJ_k)$ for all $k$, whence $\dl(\Mat_2(N))\le3\jl(N)$. The constant $3$ cannot be improved: a field $K$ of characteristic $2$ satisfies $K\circ K=0$, so $\jl(K)=1$, whereas $\dl(\Mat_2(K))=3$ by Proposition~\ref{prop:gl2}(c).

(2) By Lemma~\ref{lem:jordan-series}(b), $\jl(N)=\dlL(N)$, so (2) is a special case of (1).

(3) If $\Mat_2(N)$ is Lie solvable, then $\jl(N)<\infty$ by (1), and $N$ is nilpotent by Lemma~\ref{lem:jordan-series}(c). Conversely, if $N^t=0$, then $\Mat_2(N)^t\subseteq\Mat_2(N^t)=0$, and \ref{D:ideal} gives $\dl(\Mat_2(N))\le\lceil\log_2t\rceil$. The first inequality in the display is Lemma~\ref{lem:jordan-series}(c), the second is (1), and the last holds because the first means $t(N)\le3^{\jl(N)}$. If $N$ is commutative, then $\cJ_k(N)=N^{2^k}$ by Lemma~\ref{lem:jordan-series}(c), so $\jl(N)=\lceil\log_2t(N)\rceil$, and the display gives $\jl(N)=\dl(\Mat_2(N))=\lceil\log_2t(N)\rceil$.
\end{proof}

The following remark explains how Theorem~\ref{thm:main} relates to the classical commutative case; the sharpness of the remaining bounds is discussed in Remark~\ref{rem:sharpness}.

\begin{remark}\label{rem:two}
For algebras over fields of characteristic different from $2$, part (3) of Theorem~\ref{thm:main} extends \cite[Lemma~1.2]{PPS73} for $n=2$ from commutative to arbitrary algebras, and determines the derived length in the commutative case. More generally, Theorem~\ref{thm:main} explains the role of the prime $2$ in Theorem~\ref{thm:PPS}(a): the hypothesis on the characteristic is replaced by an intrinsic condition on $N$, the termination of its Jordan ideal series. When $2$ is invertible on $N$ this condition is nilpotency, and when $2N=0$ it is strong Lie solvability, a much weaker condition, as part (b) of the following example shows.
\end{remark}

The dichotomy between parts (2) and (3) of Theorem~\ref{thm:main} is illustrated by the following two examples of nil rings that are not nilpotent; in part (a) we also determine both invariants exactly for the Grassmann algebras on finitely many generators.

\begin{example}\label{ex:dichotomy}
(a) Let $K$ be a field of characteristic different from $2$, let $\Lambda$ be the Grassmann algebra over $K$ on generators $e_1,e_2,\dots$ subject to $e_ie_j=-e_je_i$, and let $\Lambda^+$ be the ideal generated by the $e_i$. Then $\Lambda^+$ is nil and Lie nilpotent, but not nilpotent, since $e_1e_2\cdots e_n\ne0$ for all $n$. By Theorem~\ref{thm:main}(3), $\Mat_2(\Lambda^+)$ is not Lie solvable; explicitly, $\cJ_k(\Lambda^+)$ is spanned by the monomials of degree at least $2^{k+1}-1$, since a product $u\circ v$ of monomials is $2uv$ if $u$ or $v$ has even degree and $0$ otherwise. For the Grassmann algebra $\Lambda_g$ on $g$ generators, $t(\Lambda_g^+)=g+1$, and the same description gives $\jl(\Lambda_g^+)=\lceil\log_2(g+2)\rceil-1$, while
\[
\dl(\Mat_2(\Lambda_g^+))=\lceil\log_2(g+1)\rceil .
\]
The upper bound is Theorem~\ref{thm:main}(3). For the lower bound, note that monomials of even degree are central in $\Lambda$. For monomials $u,v$ of degree $1$, the relations $[uE_{11},vE_{11}]=2uvE_{11}$, $[uE_{22},vE_{22}]=2uvE_{22}$, $[uE_{11},vE_{12}]=uvE_{12}$ and $[uE_{22},vE_{21}]=uvE_{21}$ show that $\Mat_2(\Lambda^+)^{(1)}$ contains $wE_{12}$, $wE_{21}$ and $wh$ for every monomial $w$ of degree $2$, where $h=E_{11}-E_{22}$; here we use that the terms of the derived series are $K$-subspaces. If $\Mat_2(\Lambda^+)^{(k)}$ contains these elements for every monomial $w$ of degree $2^k$, then, for monomials $u,v$ of degree $2^k$, the relations $[uE_{12},vE_{21}]=uvh$, $[uh,vE_{12}]=2uvE_{12}$ and $[uh,vE_{21}]=-2uvE_{21}$ show that $\Mat_2(\Lambda^+)^{(k+1)}$ contains them for every monomial of degree $2^{k+1}$, since every such monomial is, up to sign, a product $uv$ of this kind and $2$ is invertible in $K$. The same argument applies to $\Lambda_g$, so $\Mat_2(\Lambda_g^+)^{(k)}\ne0$ whenever $2^k\le g$, which gives the lower bound.

(b) Let $N$ be the augmentation ideal of $\F_2[x_1,x_2,\dots]/(x_1^2,x_2^2,\dots)$. Then $N$ is commutative and nil, since $y^2=0$ for every $y\in N$, but not nilpotent; nevertheless $\dl(\Mat_2(N))=3$. Indeed, $N$ is anticommutative, so $\dl(\Mat_2(N))\le3$ by Lemma~\ref{lem:anticomm}; and $x_1x_2E_{12}=[x_1E_{11},x_2E_{12}]$ and $x_3x_4E_{21}=[x_3E_{22},x_4E_{21}]$ lie in $\Mat_2(N)^{(1)}$, so that $[x_1x_2E_{12},x_3x_4E_{21}]=x_1x_2x_3x_4\cdot1\ne0$ lies in $\Mat_2(N)^{(2)}$. In characteristic $2$, therefore, nilpotency of $N$ is far from necessary.
\end{example}

We conclude this section with matrices of size at least $3$. Here the commutators $[xE_{ij},yE_{jl}]$ with pairwise distinct indices bring every power of $N$ into the derived series, and the answer, which does not depend on the characteristic, shows that the Jordan ideal series is a phenomenon of size $2$.

\begin{proposition}\label{prop:Mn}
Let $N$ be a ring and $n\ge3$. Then $\Mat_n(N)$ is Lie solvable if and only if $N$ is nilpotent, and in this case $\dl(\Mat_n(N))=\lceil\log_2t(N)\rceil$.
\end{proposition}

\begin{proof}
Let $M=\Mat_n(N)$. We show by induction on $r$ that $N^{2^r}E_{il}\subseteq M^{(r)}$ for all $r\ge0$ and all $i\ne l$. This is clear for $r=0$. Suppose that it holds for $r$, let $i\ne l$, and choose $j\notin\{i,l\}$, which is possible because $n\ge3$. For $x,y\in N^{2^r}$ the elements $xE_{ij}$ and $yE_{jl}$ lie in $M^{(r)}$, and $[xE_{ij},yE_{jl}]=xyE_{il}$ because $E_{jl}E_{ij}=0$. Since $N^{2^{r+1}}$ is the additive subgroup generated by these products $xy$, we obtain $N^{2^{r+1}}E_{il}\subseteq M^{(r+1)}$. On the other hand, \ref{D:ideal}, applied to $M$ as an ideal of itself, gives $M^{(r)}\subseteq M^{2^r}\subseteq\Mat_n(N^{2^r})$. Hence $M^{(r)}=0$ if and only if $N^{2^r}=0$, that is, if and only if $N$ is nilpotent with $t(N)\le2^r$, and both assertions follow.
\end{proof}

Thus, for every $n\ge2$, the Lie solvability of $\Mat_n(N)$ is governed by a single descending series of ideals of $N$: the Jordan ideal series when $n=2$, and the series of the powers $N^{2^r}$ when $n\ge3$. When $2$ is invertible on $N$, both conditions amount to nilpotency of $N$, and for commutative $N$ the derived lengths of $\Mat_2(N)$ and $\Mat_n(N)$ coincide by Theorem~\ref{thm:main}(3). When $2N=0$, by contrast, the two series can behave very differently: for the ring $N$ of Example~\ref{ex:dichotomy}(b), the ring $\Mat_2(N)$ has derived length $3$, whereas $\Mat_n(N)$ is not Lie solvable for $n\ge3$.

\section{Unital rings, matrix units and group algebras}\label{s:unital}

We now specialize Theorem~\ref{thm:main} to unital coefficient rings. For a unital ring $B$ of characteristic $\ell$, the subring of $\Mat_2(B)$ generated by the matrix units $E_{ij}$ is the arithmetic skeleton $\Mat_2(\Z\cdot1_B)\cong\Mat_2(\Z/\ell\Z)$, and Theorem~\ref{thm:main} combines with Proposition~\ref{prop:gl2} to separate an arithmetic condition from a structural one. Part (b) of the following corollary improves the upper bound when a nilpotent ideal with commutative quotient is available, as for triangular rings and group algebras; this part does not require an identity.

\begin{corollary}\label{cor:M2B}
Let $B$ be a ring.
\begin{enumerate}[label=\textup{(\alph*)}]
\item Suppose that $B\ne0$ is unital. Then $\Mat_2(B)$ is Lie solvable if and only if $2^mB=0$ for some $m\ge1$ and $B$ is strongly Lie solvable. In this case, if $2^m$ is the characteristic of $B$, then
\[
\max\bigl\{d_m,\ \dlL(B/2B)+1\bigr\}\le\dl(\Mat_2(B))\le d_m\cdot\dlL(B),
\]
and $\dlL(B)\le\dlL(B/2B)+\lceil\log_2m\rceil$.
\item If $B$ has an ideal $I$ with $I^t=0$ such that $B/I$ is commutative and $2^m(B/I)=0$ for some $m\ge1$, then $\dl(\Mat_2(B))\le d_m+\lceil\log_2t\rceil$.
\end{enumerate}
\end{corollary}

\begin{proof}
(a) Suppose that $\Mat_2(B)$ is Lie solvable of derived length $d$, and let $\ell$ be the characteristic of $B$. Since $\Mat_2(\Z/\ell\Z)$ is a subring of $\Mat_2(B)$, \ref{D:sub} and Proposition~\ref{prop:gl2}(b) give $\ell=2^m$ with $m\ge1$ and $d_m\le d$. Put $\bar B=B/2B$; this is a nonzero unital ring, since $B=2B$ would imply $B=2^mB=0$. By \ref{D:hom}, $\dl(\Mat_2(\bar B))\le d$, so Theorem~\ref{thm:main}(4) and Lemma~\ref{lem:jordan-series}(b) give $s:=\dlL(\bar B)=\jl(\bar B)\le d-1$. Since the projection $B\to\bar B$ maps $\Delta^{(s)}(B)$ onto $\Delta^{(s)}(\bar B)=0$, we have $\Delta^{(s)}(B)\subseteq2B$. If $\Delta^{(k)}(B)\subseteq2^aB$, then $[\Delta^{(k)}(B),\Delta^{(k)}(B)]\subseteq2^{2a}B$, so $\Delta^{(k+1)}(B)\subseteq2^{2a}B$; hence $\Delta^{(s+j)}(B)\subseteq2^{2^j}B$ for all $j$, and $\Delta^{(s+\lceil\log_2m\rceil)}(B)\subseteq2^mB=0$.

Conversely, suppose that $2^mB=0$ and $\dlL(B)=s<\infty$. For every $k$ the ring $\Delta^{(k)}(B)/\Delta^{(k+1)}(B)$ is commutative and annihilated by $2^m$, so Proposition~\ref{prop:gl2}(c) and \ref{D:hom} give $\Mat_2\bigl(\Delta^{(k)}(B)\bigr)^{(d_m)}\subseteq\Mat_2\bigl(\Delta^{(k+1)}(B)\bigr)$. By induction, $\Mat_2(B)^{(d_mk)}\subseteq\Mat_2(\Delta^{(k)}(B))$ for all $k$, whence $\dl(\Mat_2(B))\le d_ms$.

(b) By Proposition~\ref{prop:gl2}(c), $\dl(\Mat_2(B/I))\le d_m$, so \ref{D:hom}, applied to the projection $\Mat_2(B)\to\Mat_2(B/I)$, gives $\Mat_2(B)^{(d_m)}\subseteq J:=\Mat_2(I)$. This is an ideal of $\Mat_2(B)$ with $J^t\subseteq\Mat_2(I^t)=0$, so by \ref{D:sub} and \ref{D:ideal} we obtain $\Mat_2(B)^{(d_m+k)}\subseteq J^{(k)}\subseteq J^{2^k}=0$ as soon as $2^k\ge t$.
\end{proof}

For triangular rings in characteristic $2$ the derived length can be computed exactly. The following proposition shows that the bound of Corollary~\ref{cor:M2B}(b) is attained, and that the lower bound of Theorem~\ref{thm:main}(4) misses the true value by exactly one for an infinite family of rings.

\begin{proposition}\label{prop:triangular}
Let $A$ be a unital commutative ring of characteristic $2$, let $n\ge1$, and let $\Tri_n(A)$ be the ring of upper triangular $n\times n$ matrices over $A$. Then
\[
\dlL(\Tri_n(A))=1+\lceil\log_2n\rceil\qquad\text{and}\qquad\dl\bigl(\Mat_2(\Tri_n(A))\bigr)=3+\lceil\log_2n\rceil .
\]
\end{proposition}

\begin{proof}
Let $J$ be the ideal of strictly upper triangular matrices in $\Tri_n(A)$. Since $[E_{ii},aE_{ij}]=aE_{ij}$ for $i<j$, the ideal generated by $[\Tri_n(A),\Tri_n(A)]\subseteq J$ is $J$. For $a\ge1$, the set $[J^a,J^a]$ contains $[E_{i,i+a},bE_{i+a,j}]=bE_{ij}$ whenever $j-i\ge2a$ and $b\in A$, so the ideal it generates is $J^{2a}$. Hence $\Delta^{(k)}(\Tri_n(A))=J^{2^{k-1}}$ for $k\ge1$, and $\dlL(\Tri_n(A))=1+\lceil\log_2n\rceil$, because $J^r=0$ exactly when $r\ge n$.

For the second formula we identify $L=\Mat_2(\Tri_n(A))$ with the ring of block upper triangular $n\times n$ matrices with entries in $Q=\Mat_2(A)$, write $c\otimes E_{ij}$ for the block matrix with $c\in Q$ in block $(i,j)$, and let $L_0$ and $U$ be the block diagonal and the block strictly upper triangular parts of $L$; the case $n=1$ is Proposition~\ref{prop:gl2}(c), so let $n\ge2$. By Proposition~\ref{prop:gl2}(a), and since $h=E_{11}+E_{22}=1$ in characteristic $2$, the first two terms of the derived series of $Q$ are $Q^{(1)}=Ax+Ay+A\cdot1$ and $Q^{(2)}=A\cdot1$. Let $q_i=1\otimes E_{ii}$. For $i<j$ and $c\in Q$ we have $[q_i,c\otimes E_{ij}]=c\otimes E_{ij}$, and $q_i$ lies in $(Q^{(1)})^{\oplus n}$ and in $(Q^{(2)})^{\oplus n}$, the block diagonal matrices with entries in $Q^{(1)}$ and in $Q^{(2)}$. Since $[L_0,L_0]$ is computed blockwise and $[L_0,U]+[U,U]\subseteq U$, it follows that
\[
L^{(1)}=(Q^{(1)})^{\oplus n}\oplus U,\qquad L^{(2)}=(Q^{(2)})^{\oplus n}\oplus U,\qquad L^{(3)}=U .
\]
Finally, $U^{(k)}=U^{2^k}$ for all $k$: the inclusion $\subseteq$ is \ref{D:ideal}, and conversely $[U^a,U^a]$ contains the commutator of $c\otimes E_{i,i+a}$ with $1\otimes E_{i+a,j}$, which is $c\otimes E_{ij}$, for all $c\in Q$ and $j-i\ge2a$. Since $U^{2^k}=0$ exactly when $2^k\ge n$, we conclude that $\dl(L)=3+\dl(U)=3+\lceil\log_2n\rceil$.
\end{proof}

We can now describe the matrix units of a Lie solvable ring. The following theorem is obtained by applying Corollary~\ref{cor:M2B}(a) to the corner $eRe\cong\Mat_2(B)$ given by \eqref{eq:corner}; part (1) is the arithmetic content and parts (2) and (3) the structural content, and for (3) we use a classical lemma on prime rings \cite[Lemma~4]{May82}.

\begin{theorem}\label{thm:units}
Let $R$ be a Lie solvable ring of derived length $d$ containing a system of $2\times2$ matrix units $(e_{ij})$, and put $e=e_{11}+e_{22}$ and $B=e_{11}Re_{11}$.
\begin{enumerate}[label=\textup{(\arabic*)}]
\item $d\ge3$, and $e_{11}$ has additive order $2^m$ with
\[
1\le m\le2\lfloor2^{d-1}/3\rfloor,\qquad\text{equivalently}\qquad d\ge d_m .
\]
For every $d\ge3$ equality holds for $R=\Mat_2(\Z/2^m\Z)$ with $m=2\lfloor2^{d-1}/3\rfloor$.
\item The corner ring $B$ is strongly Lie solvable, with $\dlL(B/2B)\le d-1$ and $\dlL(B)\le d-1+\lceil\log_2m\rceil$.
\item For every prime ideal $P$ of $B$, the ring $B/P$ is commutative of characteristic $2$. Equivalently, $[B,B]+2B\subseteq\Nil(B)$, and $B/\Nil(B)$ is a reduced commutative ring of characteristic $2$, where $\Nil(B)$ is the prime radical of $B$.
\end{enumerate}
\end{theorem}

\begin{proof}
By \eqref{eq:corner}, $eRe\cong\Mat_2(B)$, and this is a subring of $R$, so $\dl(\Mat_2(B))\le d$ by \ref{D:sub}. The additive order of $e_{11}$ is the characteristic of the ring $B$ with identity $e_{11}$, so Corollary~\ref{cor:M2B}(a) shows that it equals $2^m$ with $m\ge1$ and $d_m\le d$, and gives (2). Since $d_m\ge3$, we get $d\ge3$. By Proposition~\ref{prop:gl2}(a), the sequence $(\alpha_k)$ is nondecreasing, so $m\le\alpha_{d_m}\le\alpha_d=2\lfloor2^{d-1}/3\rfloor$; conversely, $m\le\alpha_d$ implies $d_m\le d$ by definition. For the sharpness, $(\alpha_k)_{k\ge2}$ is strictly increasing, so for $d\ge3$ and $m=\alpha_d$ we have $\alpha_{d-1}<m\le\alpha_d$, hence $d_m=d$, and Proposition~\ref{prop:gl2}(b) applies.

For (3), let $P$ be a prime ideal of $B$. Then $B/P$ is strongly Lie solvable, because $\Delta^{(k)}(B)$ maps onto $\Delta^{(k)}(B/P)$. If $k\ge1$ is minimal with $\Delta^{(k)}(B/P)=0$, then $I=\Delta^{(k-1)}(B/P)$ is a nonzero ideal with $[I,I]=0$, and a prime ring containing a nonzero commutative right ideal is commutative \cite[Lemma~4]{May82}; hence $B/P$ is commutative. Since $2^m\cdot1=0$ in the prime ring $B/P$, whose identity has infinite or prime additive order, $B/P$ has characteristic $2$. Finally, $B/\Nil(B)$ embeds in $\prod_PB/P$, so it is commutative, semiprime, hence reduced, and of characteristic $2$.
\end{proof}

Since the prime ideals of $\Mat_2(B)$ are the ideals $\Mat_2(P)$ with $P$ prime in $B$, Theorem~\ref{thm:units} shows that every quotient of the corner $eRe$ by a prime ideal is a matrix ring $\Mat_2(D)$ over a commutative domain $D$ of characteristic $2$, in accordance with the exceptional case of the Lie theory of prime rings \cite[Theorem~4]{LM72}. The next remark compares Theorem~\ref{thm:units} with the structure theorem of \cite{ZS82} and with commutativity.

\begin{remark}\label{rem:units}
(a) That the additive order of $e_{11}$ is a power of $2$ can also be derived from the structure theorem of \cite{ZS82}, in the form \cite[Theorem~2.5]{SU17}. Put $u=[e_{12},e_{21}]=e_{11}-e_{22}$; then $[u,[e_{11},e_{12}]]=2e_{12}\in[[R,R],[R,R]]$ and $[2e_{12},e_{21}]=2u$, so the nilpotent ideal $I$ generated by $[[[R,R],[R,R]],R]$ contains $e_{11}(2u)e_{11}=2e_{11}$, and $I^t=0$ gives $2^te_{11}=0$. This yields only $m\le t$, where $t$ is the index of nilpotency of $I$, for which the available estimate is of order $10^d$ when $3$ is invertible \cite[Theorem~2.1]{SS85}; Theorem~\ref{thm:units}(1) is sharp and its proof is elementary.

(b) The corner ring in Theorem~\ref{thm:units} need not be commutative. For $R=\Mat_2(\Tri_n(\Z/2^m\Z))$ with $n\ge2$, the corner of the system $(E_{ij})$ is $\Tri_n(\Z/2^m\Z)$, and $R$ is Lie solvable by Corollary~\ref{cor:M2B}(b).
\end{remark}

We now return to Question~\ref{q:inverse}, using the notions recalled in Section~\ref{s:prelim}. Combining Theorem~\ref{thm:units} with \ref{D:units}, we obtain the following answer; the matrix units in part (a) are constructed from an element $s$ and its inverse, whose idempotents $\dd(s)$ and $\rr(s)$ commute but differ.

\begin{corollary}\label{cor:inverse}
Let $R$ be a ring.
\begin{enumerate}[label=\textup{(\alph*)}]
\item The semigroup $(R,\cdot)$ contains a non-Clifford inverse subsemigroup if and only if $R$ contains a system of $2\times2$ matrix units.
\item If $R$ is Lie solvable, then no subsemigroup of $(R,\cdot)$ is isomorphic to a Brandt semigroup $B_n$ with $n\ge3$.
\item If $R$ is Lie solvable and either has no $2$-torsion or satisfies $\dl(R)\le2$, then every inverse subsemigroup of $(R,\cdot)$ is a Clifford semigroup. Both conditions are sharp: for every field $K$ of characteristic $2$, the algebra $\Mat_2(K)$ is Lie solvable of derived length $3$ and contains the rook monoid.
\end{enumerate}
\end{corollary}

\begin{proof}
(a) Let $S$ be a non-Clifford inverse subsemigroup of $(R,\cdot)$, choose $s\in S$ with $\dd(s)\ne\rr(s)$, and put $s'=s^{-1}$, $e=s's$ and $f=ss'$; these idempotents of $S$ commute. Interchanging $s$ and $s'$, which interchanges $e$ and $f$, if necessary, we may assume that $p=e-ef\ne0$; indeed, if $e=ef$ and $f=fe$, then $e=f$. Since $e$ and $f$ are commuting idempotents, $p$ is an idempotent with $pe=ep=p$ and $pf=fp=0$. Put $e_{22}=p$, $e_{12}=sp$, $e_{21}=ps'$ and $e_{11}=sps'$. Then $e_{21}e_{12}=ps'sp=pep=p=e_{22}$ and $e_{12}e_{21}=sps'=e_{11}$, and $e_{11}$ is an idempotent because $sps'sps'=speps'=sps'$. Since $fs=s$ and $s'f=s'$, we have $fe_{11}=e_{11}=e_{11}f$, whence $e_{22}e_{11}=pfe_{11}=0$ and $e_{11}e_{22}=e_{11}fp=0$. The remaining relations follow from $e_{12}=e_{11}e_{12}=e_{12}e_{22}$ and $e_{21}=e_{22}e_{21}=e_{21}e_{11}$; for instance, $e_{12}e_{12}=e_{12}e_{22}e_{11}e_{12}=0$. Finally, $e_{22}=p\ne0$, and $p=e_{21}e_{11}e_{12}$ shows that $e_{11},e_{12},e_{21}$ are nonzero. Conversely, a system of matrix units spans, together with $0$, a copy of $B_2$.

(b) If $\{f_{ij}\}\cup\{z\}\cong B_n$ with zero $z$, then $z$ is an idempotent with $zf_{ij}=f_{ij}z=z$, and the nonzero elements $f_{ij}-z$ satisfy $(f_{ij}-z)(f_{kl}-z)=\delta_{jk}(f_{il}-z)$; they form a system of $n\times n$ matrix units, which is impossible for $n\ge3$ by \ref{D:units}.

(c) By (a) and Theorem~\ref{thm:units}(1), a non-Clifford inverse subsemigroup forces $\dl(R)\ge3$ and $2^me_{11}=0$ for some $m\ge1$ and some nonzero idempotent $e_{11}$, which is impossible under either condition. The example is Proposition~\ref{prop:gl2}(c).
\end{proof}

We turn to group algebras. In characteristic $2$, Theorem~\ref{thm:main}(2) reduces the Lie solvability of $\Mat_2(N)$ to the strong Lie solvability of $N$, which is well understood for group algebras; combined with Theorem~\ref{thm:groupalg}, it gives the following purely group-theoretic criterion.

\begin{corollary}\label{cor:M2KG}
Let $K$ be a field of characteristic $2$ and $G$ a group. Then $\Mat_2(KG)$ is Lie solvable if and only if $G'$ is a finite $2$-group, and in this case
\[
\max\bigl\{3,\ \lceil\log_2(t(G')+1)\rceil+1\bigr\}\le\dl(\Mat_2(KG))\le\lceil\log_2t(G')\rceil+3 .
\]
\end{corollary}

\begin{proof}
Since $2KG=0$, Theorem~\ref{thm:main}(2) shows that $\Mat_2(KG)$ is Lie solvable if and only if $KG$ is strongly Lie solvable, that is, by Theorem~\ref{thm:groupalg}(b), if and only if $G'$ is a finite $2$-group. The lower bound $3$ is Corollary~\ref{cor:M2B}(a) with $m=1$, since $d_1=3$. If $G$ is not abelian, the other lower bound follows from Theorem~\ref{thm:main}(4), since $\jl(KG)=\dlL(KG)$ by Lemma~\ref{lem:jordan-series}(b), and Theorem~\ref{thm:groupalg}(c); if $G$ is abelian, then $t(G')=1$, and the two bounds coincide, so that $\dl(\Mat_2(KG))=3$. For the upper bound, apply Corollary~\ref{cor:M2B}(b) with $m=1$ to the ideal $I=\omega(KG')KG$: the quotient $KG/I\cong K[G/G']$ is commutative of characteristic $2$, and $I^{t(G')}=\omega(KG')^{t(G')}KG=0$.
\end{proof}

Corollary~\ref{cor:M2KG} should be compared with Theorem~\ref{thm:PPS}(b), according to which $KG$ is Lie solvable as soon as $G$ has a $2$-abelian subgroup of index at most $2$. For instance, $KS_3$ and the algebras $KD_{2^n}$, where $D_{2^n}$ is the dihedral group of order $2^n$, are Lie solvable of derived length at most $3$, by \cite[Proposition~1]{Ju06b} applied to their cyclic subgroups of index $2$. However, $\Mat_2(KS_3)$ is not Lie solvable, since $S_3'=A_3$ has order $3$, while for $n\ge3$ the group $D_{2^n}'$ is cyclic of order $2^{n-2}$, so that $t(D_{2^n}')=2^{n-2}$ and $n\le\dl(\Mat_2(KD_{2^n}))\le n+1$. Thus the derived lengths of the algebras $KD_{2^n}$ are bounded, whereas those of the algebras $\Mat_2(KD_{2^n})$ are not.

Both bounds in Corollary~\ref{cor:M2KG} are attained. Since the derived series is compatible with extension of scalars, $\dl(\Mat_2(KG))$ does not depend on the field $K$ of characteristic $2$, and a direct computation over $\F_2$ shows that it equals the upper bound, namely $4$, $4$, $5$ and $5$, for $G=D_8$, the quaternion group $Q_8$, the alternating group $A_4$ and $D_{16}$, respectively, and the lower bound $5$ for $G=D_{32}$. The computation is exact linear algebra over $\F_2$: each term of the derived series of $\Mat_2(\F_2G)$ is spanned by the commutators of pairs of basis elements of the previous term, and a row-reduced basis of it is computed until the zero space is reached. In particular, each of the inequalities $n\le\dl(\Mat_2(KD_{2^n}))\le n+1$ is an equality for some $n$.

Since $\Mat_2(A)\cong\Mat_2(K)\otimes_KA$ for every $K$-algebra $A$, Theorem~\ref{thm:main}(2) completely settles the question studied in [30] when one factor is
Lie solvable if and only if $A$ is strongly Lie solvable. The following example makes the failure of Lie solvability for tensor products concrete.

\begin{example}\label{ex:tensor}
Let $K$ be a field of characteristic $2$. The homomorphism $KS_3\to\Mat_2(K)$ induced by the isomorphism $S_3\cong\mathrm{GL}_2(\F_2)$ is surjective, since the six matrices of $\mathrm{GL}_2(\F_2)$ span $\Mat_2(\F_2)$. Hence $\dl(KS_3)\ge\dl(\Mat_2(K))=3$ by \ref{D:hom} and Proposition~\ref{prop:gl2}(c), and the reverse inequality was recalled above, so $KS_3$ is Lie solvable of derived length $3$; for $K=\F_2$ this is also noted in \cite[p.~291]{Sp08}. However, $KS_3$ is not strongly Lie solvable, so $\Mat_2(K)\otimes_KKS_3\cong\Mat_2(KS_3)$ is not Lie solvable, although both factors are Lie solvable of derived length $3$. This can also be seen directly: the surjection above induces a surjection of $\Mat_2(KS_3)$ onto $\Mat_2(\Mat_2(K))\cong\Mat_4(K)$, which is not Lie solvable by \ref{D:units}. The examples in \cite{Ri96} show that the main theorem there is the best possible of its kind; the present example shows, in the same direction, that the hypothesis that each factor be Lie nilpotent or Lie metabelian cannot be weakened to derived length $3$, even when one factor is a group algebra.
\end{example}

\section{Steinberg algebras, Leavitt path algebras and inverse semigroup algebras}\label{s:steinberg}

In this section we apply the main theorem to the Steinberg algebras of ample groupoids recalled in Section~\ref{s:prelim}. Orbits with two points produce blocks $\Mat_2(KH)$, to which Theorem~\ref{thm:main} and Corollary~\ref{cor:M2KG} apply, while larger orbits are excluded by \ref{D:units}. The first part of the following characterization, valid in every characteristic, rests on the observation that $n$ distinct points of one orbit produce a system of $n\times n$ matrix units in $A_K(\cG)$; the remaining parts combine Proposition~\ref{prop:finite-orbits} with Sections~\ref{s:main} and~\ref{s:unital}.

\begin{theorem}\label{thm:steinberg}
Let $\cG$ be an ample groupoid and $K$ a field.
\begin{enumerate}[label=\textup{(\arabic*)}]
\item If $A_K(\cG)$ is Lie solvable, then every orbit of $\cG$ has at most two points; if, moreover, $\chr K\ne2$, then $\cG$ is a group bundle.
\item If $\chr K\ne2$, then $A_K(\cG)$ is Lie solvable if and only if $\cG$ is a group bundle and $\sup_x\dl(K[\cG_x^x])<\infty$, and then $\dl(A_K(\cG))=\sup_x\dl(K[\cG_x^x])$. If $\chr K=0$, this happens if and only if $A_K(\cG)$ is commutative.
\item If $\chr K=2$, then $A_K(\cG)$ is Lie solvable if and only if the following conditions hold:
\begin{enumerate}[label=\textup{(\roman*)}]
\item every orbit has at most two points;
\item the commutator subgroups $H_O'$ of the isotropy groups at the two-point orbits $O$ are finite $2$-groups of uniformly bounded order;
\item $\sup_x\dl(K[\cG_x^x])<\infty$, the supremum being taken over the points whose orbit is a single point.
\end{enumerate}
\end{enumerate}
\end{theorem}

\begin{proof}
(1) Let $x_1,\dots,x_n$ be $n\ge2$ distinct points in one orbit. Choose $\gamma_j\in\cG$ with $s(\gamma_j)=x_1$ and $r(\gamma_j)=x_j$, where $\gamma_1=x_1$, and pairwise disjoint compact open sets $V_1,\dots,V_n\subseteq\cG^{(0)}$ with $x_j\in V_j$. Choose compact open bisections $U_j\ni\gamma_j$ with $U_1=V_1$, and put $U'_j=V_jU_jV_1$, a compact open bisection containing $\gamma_j$ with $s(U'_j)\subseteq V_1$ and $r(U'_j)\subseteq V_j$. The set $W=\bigcap_js(U'_j)$ is compact open, being a finite intersection of compact open subsets of a Hausdorff space, and contains $x_1$. Put $W_j=U'_jW$, a compact open bisection containing $\gamma_j$ with $s(W_j)=W$ and $r(W_j)\subseteq V_j$, and let $a_{ij}=1_{W_iW_j^{-1}}$. Since $W_j^{-1}W_j=W$ and $W_iW=W_i$, we have $W_iW_j^{-1}W_jW_l^{-1}=W_iW_l^{-1}$, while $W_j^{-1}W_k=\emptyset$ for $j\ne k$, because $r(W_j)\subseteq V_j$ and $r(W_k)\subseteq V_k$ are disjoint. Hence $a_{ij}a_{kl}=\delta_{jk}a_{il}$, and $a_{ij}\ne0$ since $\gamma_i\gamma_j^{-1}\in W_iW_j^{-1}$; thus the $a_{ij}$ form a system of $n\times n$ matrix units in $A_K(\cG)$. An orbit with at least three points is therefore excluded by \ref{D:units}. If $\chr K\ne2$, an orbit with two points yields a system of $2\times2$ matrix units in a $K$-algebra, in which the additive order of $e_{11}$ is infinite or an odd prime; this contradicts Theorem~\ref{thm:units}(1).

(2) If $A_K(\cG)$ is Lie solvable, then $\cG$ is a group bundle by (1). Suppose now that $\cG$ is a group bundle. Then every orbit is a single point $\{x\}$, and $\cG|_{\{x\}}=\cG_x^x$, so Proposition~\ref{prop:finite-orbits} gives $\dl(A_K(\cG))=\sup_x\dl(K[\cG_x^x])$. Hence $A_K(\cG)$ is Lie solvable if and only if this supremum is finite, which proves both implications and the formula. If $\chr K=0$ and $A_K(\cG)$ is Lie solvable, then every $K[\cG_x^x]$ is Lie solvable, so every $\cG_x^x$ is abelian by Theorem~\ref{thm:PPS}(b), and $A_K(\cG)$ is a subdirect product of commutative algebras, hence commutative; the converse is clear.

(3) If $A_K(\cG)$ is Lie solvable, then (i) holds by (1), all orbits are finite, and $\dl(\Mat_{|O|}(KH_O))\le\dl(A_K(\cG))$ for every orbit $O$ by Proposition~\ref{prop:finite-orbits}; this gives (iii). For a two-point orbit $O$, Corollary~\ref{cor:M2KG} shows that $H_O'$ is a finite $2$-group with $\log_2t(H_O')<\dl(A_K(\cG))$, and $t(H)\ge1+\log_2|H|$ for every finite $2$-group $H$ (Theorem~\ref{thm:groupalg}(a)) bounds $|H_O'|$ uniformly, which is (ii). Conversely, under (i)--(iii), Corollary~\ref{cor:M2KG} and $t(H)\le|H|$ give $\dl(\Mat_2(KH_O))\le\lceil\log_2|H_O'|\rceil+3$ for the two-point orbits, and Proposition~\ref{prop:finite-orbits} shows that $A_K(\cG)$ is Lie solvable.
\end{proof}

We now make Theorem~\ref{thm:steinberg} explicit for Leavitt path algebras. Let $E$ be a directed graph, with source and range maps $s,r$ on its edges. A vertex that emits infinitely many edges is an \emph{infinite emitter}, and $E$ is \emph{row-finite} if it has none; a vertex is \emph{singular} if it emits no edges or is an infinite emitter. A \emph{cycle} is a path $e_1\cdots e_n$ with $r(e_n)=s(e_1)$ and $s(e_i)\ne s(e_j)$ for $i\ne j$, and an \emph{exit} of this cycle is an edge $f$ with $s(f)=s(e_i)$ and $f\ne e_i$ for some $i$. The boundary path space $\partial E$ consists of the infinite paths and of the finite paths, vertices included, that end at a singular vertex. For $x\in\partial E$ of length at least $n$, let $\sigma^n(x)$ be the path obtained by deleting the first $n$ edges. The boundary path groupoid $\cG_E=\{(x,m-n,y):\sigma^m(x)=\sigma^n(y)\}$, with $s(x,k,y)=y$ and $r(x,k,y)=x$, is ample, its orbits are the tail-equivalence classes of $\partial E$, where $x,y\in\partial E$ are \emph{tail equivalent} if $\sigma^m(x)=\sigma^n(y)$ for some $m,n\ge0$, its isotropy groups are subgroups of $\Z$, and the Leavitt path algebra $L_K(E)$, for whose definition we refer to \cite[Definition~1.2.3]{AAS17}, is isomorphic to $A_K(\cG_E)$ by \cite[Example~3.2]{CS15} (see \cite[Proposition~4.3]{CFST14} for an earlier special case); this isomorphism is all that we use about $L_K(E)$. Since Lie solvability of a given derived length is a polynomial identity, Lie solvable Leavitt path algebras are PI-algebras, which were characterized in \cite[Theorem~3.1]{BLR16}. It is proved in \cite{NZ22} that, if $\chr K\ne2$, then $L_K(E)$ is Lie solvable if and only if each connected component of $E$ is an isolated vertex or a vertex carrying exactly one loop and no other edge, in which case $L_K(E)$ is commutative; the case of characteristic $2$ is described there by a combinatorial condition on graphs without exits, and the solvable index of Lie solvable Leavitt path algebras is computed. For the Lie algebra of skew-symmetric elements and the Jordan algebra of symmetric elements with respect to the standard involution, the corresponding questions were answered in \cite[Theorems~4.6 and~4.11]{KD25}. 
In characteristic different from $2$, the following corollary gives a groupoid proof of the nontrivial implication of the characterization of [28], whose only graph-theoretic input is the determination of the graphs whose boundary path groupoid is a group bundle. In characteristic $2$, it gives a groupoid description of Lie solvability together with a short computation of the derived length.

\begin{corollary}\label{cor:LPA}
Let $K$ be a field and $E$ a directed graph.
\begin{enumerate}[label=\textup{(\alph*)}]
\item If $\chr K\ne2$ and $L_K(E)$ is Lie solvable, then each connected component of $E$ is an isolated vertex or a vertex carrying exactly one loop and no other edge, and $L_K(E)$ is a direct sum of copies of $K$ and of $K[x,x^{-1}]$, hence commutative.
\item If $\chr K=2$, then $L_K(E)$ is Lie solvable if and only if every tail-equivalence class of $\partial E$ has at most two elements. In this case $\dl(L_K(E))=3$ if some class has exactly two elements, and $L_K(E)$ is commutative otherwise.
\end{enumerate}
\end{corollary}

\begin{proof}
(a) By Theorem~\ref{thm:steinberg}(1), $\cG_E$ is a group bundle. If some edge $e$ satisfies $s(e)\ne r(e)$, choose $\mu\in\partial E$ starting at $r(e)$ (take $\mu=r(e)$ if $r(e)$ is singular, and otherwise choose edges successively until a singular vertex is reached or an infinite path is obtained). Then $e\mu\in\partial E$ and $\sigma^1(e\mu)=\mu\ne e\mu$, so the orbit of $\mu$ has two points, which is impossible. Hence every edge is a loop. If a vertex emits two loops $e\ne f$, then $x=eee\cdots$ and $y=fx$ are distinct and tail equivalent, which is again impossible. Thus each connected component of $E$ is an isolated vertex or a vertex carrying exactly one loop and no other edge. Hence $L_K(E)$ is the direct sum of the Leavitt path algebras of these components \cite[Proposition~1.2.14]{AAS17}, which are $K$ for an isolated vertex and $K[x,x^{-1}]$ for a vertex with one loop \cite[Proposition~1.3.4]{AAS17}.

(b) The isotropy groups of $\cG_E$ are abelian, so conditions (ii) and (iii) of Theorem~\ref{thm:steinberg}(3) hold automatically, and the first assertion follows. Every algebra $KH_O$ is then commutative and unital of characteristic $2$, so $\dl(\Mat_2(KH_O))=3$ for every two-point orbit $O$ by Proposition~\ref{prop:gl2}(c), whereas $\dl(KH_O)\le1$ for every one-point orbit $O$. Proposition~\ref{prop:finite-orbits} therefore gives $\dl(L_K(E))=3$ if some orbit has two points; otherwise $L_K(E)$ is a subdirect product of the commutative algebras $KH_O$, hence commutative.
\end{proof}

The hypothesis $\chr K\ne2$ in Corollary~\ref{cor:LPA}(a) cannot be dropped: for the graph $A_2$ with one edge $e$ between two distinct vertices, $\partial A_2=\{r(e),e\}$ is a single orbit with two points, and $L_K(A_2)\cong\Mat_2(K)$ is Lie solvable when $\chr K=2$. We conclude with inverse semigroup algebras. For an inverse semigroup $S$, the algebra $KS$ is the Steinberg algebra of the universal groupoid of $S$ \cite[Theorem~5.3]{St10}, and for finite $S$ this yields the classical decomposition
\[
KS\cong\prod_D\Mat_{n_D}(KG_D),
\]
where $D$ runs over the $\cD$-classes of $S$, $n_D$ is the number of idempotents in $D$, and $G_D$ is a maximal subgroup contained in $D$ \cite[Theorem~4.3]{St06} (see also \cite[Corollary~5.4]{St10}). Combined with Sections~\ref{s:main} and~\ref{s:unital}, this decomposition gives the following description.

\begin{corollary}\label{cor:finite-inverse}
Let $S$ be a finite inverse semigroup and $K$ a field of characteristic $p\ge0$, with $n_D$ and $G_D$ as above.
\begin{enumerate}[label=\textup{(\alph*)}]
\item If $p\ne2$, then $KS$ is Lie solvable if and only if $S$ is a Clifford semigroup and every $G_D$ is $p$-abelian.
\item If $p=2$, then $KS$ is Lie solvable if and only if, for every $\cD$-class $D$, either $n_D=1$ and $G_D$ has a $2$-abelian subgroup of index at most $2$, or $n_D=2$ and $G_D'$ is a $2$-group.
\end{enumerate}
In both cases $\dl(KS)=\max_D\dl(\Mat_{n_D}(KG_D))$.
\end{corollary}

\begin{proof}
A finite product is Lie solvable if and only if each factor is, with derived length the maximum. If $n_D\ge3$, the factor $\Mat_{n_D}(KG_D)$ is not Lie solvable by \ref{D:units}. If $p\ne2$ and $n_D=2$, it is not Lie solvable by Corollary~\ref{cor:M2B}(a), since $2$ is invertible in $KG_D$; hence $n_D=1$ for all $D$, which means that $S$ is a Clifford semigroup, and the factors are the group algebras $KG_D$, to which Theorem~\ref{thm:PPS}(b) applies. If $p=2$, the case $n_D=2$ is Corollary~\ref{cor:M2KG} and the case $n_D=1$ is Theorem~\ref{thm:PPS}(b).
\end{proof}

For instance, if $\chr K=2$, the Brandt semigroup $B(S_3,2)$, which has a $\cD$-class with two idempotents and maximal subgroup $S_3$, has a non-Lie-solvable algebra $KB(S_3,2)\cong K\times\Mat_2(KS_3)$, although $KS_3$ is Lie solvable.

\section{Concluding remarks}\label{s:remarks}

We close with remarks on the relation between our results and the literature and on the sharpness of Theorem~\ref{thm:main}. The first remark separates what is known from what is new.

\begin{remark}[Relation with known results]\label{rem:novelty}
For commutative algebras over fields, Theorem~\ref{thm:main} recovers \cite[Lemma~1.2]{PPS73} for $n=2$. The case $m=1$ of Proposition~\ref{prop:gl2}(c) is contained in the proof of that lemma, whose computation Lemma~\ref{lem:anticomm} extends to anticommutative rings, while Proposition~\ref{prop:gl2}(c) for $m\ge2$ and Theorem~\ref{thm:main}(3) add the exact derived length. Property \ref{D:units} is classical, and Proposition~\ref{prop:finite-orbits} reduces to \cite[Proposition~3.1]{St18} when the unit space is finite. For rings without $2$-torsion, Corollary~\ref{cor:inverse}(c) can also be derived from the structure theorem of \cite{ZS82} (Remark~\ref{rem:units}(a)), and Corollary~\ref{cor:LPA}(a) is the nontrivial implication of the characterization of \cite{NZ22} in characteristic different from $2$. The solvable index of Lie solvable Leavitt path algebras is also computed in \cite{NZ22}, and for row-finite graphs the value of the derived length in Corollary~\ref{cor:LPA}(b) can be read off from the structure theorem in characteristic $2$ given there; what our contribution there is a groupoid proof, which treats arbitrary graphs
without the structure theory of Leavitt path algebras. Proposition~\ref{prop:Mn}, which rests on the same identity as \ref{D:units}, Corollary~\ref{cor:inverse}(a) and the description of the graphs whose boundary path groupoid is a group bundle, used in the proof of Corollary~\ref{cor:LPA}(a), are elementary observations that may well be known; we include their proofs for completeness. To the best of our knowledge, the following results are new: Theorem~\ref{thm:main}, with the optimal constant $3$, and the Jordan ideal series on which it rests; the closed form of the derived series of $\gl_2(\Z)$ and the derived lengths of the rings $\Mat_2(\Z/2^m\Z)$ in Proposition~\ref{prop:gl2}; Corollary~\ref{cor:M2B}, Proposition~\ref{prop:triangular}, Theorem~\ref{thm:units} and Corollary~\ref{cor:M2KG}; Proposition~\ref{prop:finite-orbits} for infinite unit spaces and Theorem~\ref{thm:steinberg}; the resulting proof of Corollary~\ref{cor:LPA}(a); and Corollary~\ref{cor:finite-inverse}.
\end{remark}

The bounds of Theorem~\ref{thm:main} are of three kinds, and the examples obtained above show how far each of them is attained.

\begin{remark}[Sharpness]\label{rem:sharpness}
(a) The constant $3$ in Theorem~\ref{thm:main}(1) is best possible, since it is attained by every field $K$ of characteristic $2$, for which $\jl(K)=\dlL(K)=1$ and $\dl(\Mat_2(K))=3$. When $2N=0$ and $N\ne0$ is unital, Theorem~\ref{thm:main}(2),(4) gives $\dlL(N)+1\le\dl(\Mat_2(N))\le3\dlL(N)$. For $N=\Tri_n(A)$ with $A$ commutative of characteristic $2$, Proposition~\ref{prop:triangular} gives $\dl(\Mat_2(N))=\dlL(N)+2$, so the lower bound is missed by exactly one and the upper bound is far from the truth when $n$ is large.

(b) For rings of characteristic $2^m$ with $m\ge2$, the lower bound in Theorem~\ref{thm:main}(4) can be missed already in the commutative case: for $N=\Z/8\Z$ we have $\cJ_1(N)=2N$ and $\cJ_2(N)=0$, so $\jl(N)=2$, whereas $\dl(\Mat_2(N))=d_3=4$ by Proposition~\ref{prop:gl2}(b).

(c) When $2$ is invertible on $N$, the upper bound $\lceil\log_2t(N)\rceil$ in Theorem~\ref{thm:main}(3) is attained for commutative $N$ (Theorem~\ref{thm:main}(3)), for the Grassmann rings $\Lambda_g^+$ (Example~\ref{ex:dichotomy}(a)), and for the algebra $N$ of strictly upper triangular $n\times n$ matrices over a field of characteristic different from $2$: then $t(N)=n$, and $\Mat_2(N)$ contains the subring $\{xE_{11}:x\in N\}\cong N$, whose derived length is $\lceil\log_2n\rceil$ because $[N^a,N^a]=N^{2a}$ for all $a\ge1$. In all these cases the derived length of $\Mat_2(N)$ is determined by the index of nilpotency of $N$ alone, as it is for $\Mat_n(N)$ with $n\ge3$ by Proposition~\ref{prop:Mn}. For $N=\Lambda_2^+$ we have $t(N)=3$, $\jl(N)=1$ and $\dl(\Mat_2(N))=2$ by Example~\ref{ex:dichotomy}(a), so that every inequality in Theorem~\ref{thm:main}(3) other than $\jl(N)\le\dl(\Mat_2(N))$ is an equality. Since $\log_23<2$, the last bound in Theorem~\ref{thm:main}(3) improves the constant $3$ of Theorem~\ref{thm:main}(1) to $\log_23$, up to rounding, when $2$ is invertible on $N$, and this bound is attained when $\jl(N)=1$.
\end{remark}

\section*{Declarations}

\begin{itemize}
\item {\bf Competing interests:} The authors declare that they have no competing interests.
\item {\bf Authors' contributions:} All authors contributed equally to this work. The final editing was done by the corresponding author Nguyen Thi Mong Tuyen, and all authors approved the final manuscript.
\item {\bf Data availability:} Data sharing is not applicable to this article, as no datasets were generated or analyzed during the current study.
\item {\bf Declaration of generative AI use:} During the preparation of this work, the authors used AI-assisted tools to assist with literature searches, the checking of the mathematical arguments, and the editing of the manuscript. All suggestions incorporated into the manuscript were independently checked by the authors. The authors verified the mathematical results and proofs and take full responsibility for the final manuscript.
\end{itemize}

\section*{Data availability statement}
\addcontentsline{toc}{section}{Data availability statement}
The computations reported after Corollary~\ref{cor:M2KG} were carried out with code that is available from the corresponding author on request.


\section*{Authors' addresses}

{\small
\begin{itemize}[leftmargin=1.1em,labelwidth=0.8em,labelsep=0.3em,itemsep=1pt,topsep=0pt,parsep=0pt]
\item[\textsuperscript{a}] Faculty of Mathematics and Computer Science, University of Science, Ho Chi Minh City, Vietnam
\item[\textsuperscript{b}] Vietnam National University, Ho Chi Minh City, Vietnam
\item[\textsuperscript{c}] Faculty of Mathematics and Informatics Education, School of Education, Dong Thap University, Dong Thap Province, Vietnam
\item[\textsuperscript{d}] University of Economics and Law, Ho Chi Minh City, Vietnam
\end{itemize}

\medskip\noindent
Email: \href{mailto:mhbien@hcmus.edu.vn}{mhbien@hcmus.edu.vn} (Mai Hoang Bien), \href{mailto:ntmtuyen@dthu.edu.vn}{ntmtuyen@dthu.edu.vn} (Nguyen Thi Mong Tuyen), \href{mailto:vula@uel.edu.vn}{vula@uel.edu.vn} (Le Anh Vu)\\
ORCID: 0000-0003-4377-3786 (Mai Hoang Bien), ORCID: 0000-0002-6708-8866 (Le Anh Vu)\par}


\begin{thebibliography}{99}\small

\bibitem{AAS17} G.~Abrams, P.~Ara, M.~Siles Molina, \emph{Leavitt Path Algebras}, Lecture Notes in Math. \textbf{2191}, Springer, London, 2017.

\bibitem{AM12} G.~Abrams, Z.~Mesyan, Simple Lie algebras arising from Leavitt path algebras, J. Pure Appl. Algebra \textbf{216} (2012), 2302--2313.

\bibitem{BMM96} K.~I.~Beidar, W.~S.~Martindale III, A.~V.~Mikhalev, \emph{Rings with Generalized Identities}, Monogr. Textbooks Pure Appl. Math. \textbf{196}, Marcel Dekker, New York, 1996.

\bibitem{BLR16} J.~P.~Bell, T.~H.~Lenagan, K.~M.~Rangaswamy, Leavitt path algebras satisfying a polynomial identity, J. Algebra Appl. \textbf{15} (2016), no.~5, Article 1650084, 13~pp.

\bibitem{CS07} F.~Catino, E.~Spinelli, A note on strong Lie derived length of group algebras, Boll. Unione Mat. Ital. Sez. B (8) \textbf{10} (2007), 83--86.

\bibitem{CEPSS19} L.~O.~Clark, R.~Exel, E.~Pardo, A.~Sims, C.~Starling, Simplicity of algebras associated to non-Hausdorff groupoids, Trans. Amer. Math. Soc. \textbf{372} (2019), 3669--3712.

\bibitem{CFST14} L.~O.~Clark, C.~Farthing, A.~Sims, M.~Tomforde, A groupoid generalisation of Leavitt path algebras, Semigroup Forum \textbf{89} (2014), 501--517.

\bibitem{CH20} L.~O.~Clark, R.~Hazrat, \'Etale groupoids and Steinberg algebras, a concise introduction, in: \emph{Leavitt Path Algebras and Classical $K$-Theory}, Indian Stat. Inst. Ser., Springer, Singapore, 2020, 73--101.

\bibitem{CS15} L.~O.~Clark, A.~Sims, Equivalent groupoids have Morita equivalent Steinberg algebras, J. Pure Appl. Algebra \textbf{219} (2015), 2062--2075.


\bibitem{Her69} I.~N.~Herstein, \emph{Topics in Ring Theory}, University of Chicago Press, Chicago, 1969.

\bibitem{How95} J.~M.~Howie, \emph{Fundamentals of Semigroup Theory}, London Math. Soc. Monogr. (N.S.) \textbf{12}, Oxford University Press, Oxford, 1995.

\bibitem{Je41} S.~A.~Jennings, The structure of the group ring of a $p$-group over a modular field, Trans. Amer. Math. Soc. \textbf{50} (1941), 175--185.

\bibitem{Je42} S.~A.~Jennings, Central chains of ideals in an associative ring, Duke Math. J. \textbf{9} (1942), 341--355.

\bibitem{Je47} S.~A.~Jennings, On rings whose associated Lie rings are nilpotent, Bull. Amer. Math. Soc. \textbf{53} (1947), 593--597.

\bibitem{Ju06} T.~Juh\'asz, On the derived length of Lie solvable group algebras, Publ. Math. Debrecen \textbf{68} (2006), 243--256.

\bibitem{Ju06b} T.~Juh\'asz, Recent results on the derived length of Lie solvable group algebras, Acta Math. Acad. Paedagog. Ny\'{\i}regyh\'azi. (N.S.) \textbf{22} (2006), 237--245.

\bibitem{KD25} H.~V.~Khanh, L.~Q.~Danh, A note on Lie and Jordan structures of Leavitt path algebras, J. Algebra Appl., Article 2750156 (published online 25 February 2026), doi:10.1142/S0219498827501568.

\bibitem{KS24} M.~Khrypchenko, S.~Siciliano, Identities and derived lengths of finitary incidence algebras and their group of units, Comm. Algebra \textbf{52} (2024), 916--922.

\bibitem{Lam99} T.~Y.~Lam, \emph{Lectures on Modules and Rings}, Grad. Texts in Math. \textbf{189}, Springer, New York, 1999.

\bibitem{LM72} C.~Lanski, S.~Montgomery, Lie structure of prime rings of characteristic $2$, Pacific J. Math. \textbf{42} (1972), 117--136.

\bibitem{Loo75} O.~Loos, \emph{Jordan Pairs}, Lecture Notes in Math. \textbf{460}, Springer, Berlin, 1975.

\bibitem{May82} J.~H.~Mayne, Ideals and centralizing mappings in prime rings, Proc. Amer. Math. Soc. \textbf{86} (1982), 211--212.

\bibitem{McC04} K.~McCrimmon, \emph{A Taste of Jordan Algebras}, Universitext, Springer, New York, 2004.

\bibitem{Me13} Z.~Mesyan, Commutator Leavitt path algebras, Algebr. Represent. Theory \textbf{16} (2013), 1207--1232.

\bibitem{MSW12} J.~Meyer, J.~Szigeti, L.~van Wyk, A Cayley--Hamilton trace identity for $2\times2$ matrices over Lie-solvable rings, Linear Algebra Appl. \textbf{436} (2012), 2578--2582.

\bibitem{Na20} T.~G.~Nam, Simple Lie algebras arising from Steinberg algebras of Hausdorff ample groupoids, J. Algebra \textbf{595} (2022), 194--215.

\bibitem{NZ22} T.~G.~Nam, Z.~Zhang, Lie solvable Leavitt path algebras, J. Algebra Appl. \textbf{21} (2022), no.~10, Article 2250203.

\bibitem{PPS73} I.~B.~S.~Passi, D.~S.~Passman, S.~K.~Sehgal, Lie solvable group rings, Canad. J. Math. \textbf{25} (1973), 748--757.

\bibitem{Ri96} D.~M.~Riley, The tensor product of Lie soluble algebras, Arch. Math. (Basel) \textbf{66} (1996), 372--377.

\bibitem{Seh78} S.~K.~Sehgal, \emph{Topics in Group Rings}, Monogr. Textbooks Pure Appl. Math. \textbf{50}, Marcel Dekker, New York, 1978.

\bibitem{Sha92} A.~Shalev, The derived length of Lie soluble group rings I, J. Pure Appl. Algebra \textbf{78} (1992), 291--300.

\bibitem{SS85} R.~K.~Sharma, J.~B.~Srivastava, Lie solvable rings, Proc. Amer. Math. Soc. \textbf{94} (1985), 1--8.

\bibitem{Sh74} L.~B.~Shneperman, Maximal inverse subsemigroups of the semigroup of linear transformations, Izv. Vyssh. Uchebn. Zaved. Mat. 1974, no.~11, 93--100 (in Russian).

\bibitem{SU17} S.~Siciliano, H.~Usefi, Lie structure of smash products, Israel J. Math. \textbf{217} (2017), 93--110.

\bibitem{So02} L.~Solomon, Representations of the rook monoid, J. Algebra \textbf{256} (2002), 309--342.

\bibitem{Sp08} E.~Spinelli, Group algebras with minimal strong Lie derived length, Canad. Math. Bull. \textbf{51} (2008), 291--297.

\bibitem{St06} B.~Steinberg, M\"obius functions and semigroup representation theory, J. Combin. Theory Ser. A \textbf{113} (2006), 866--881.

\bibitem{St10} B.~Steinberg, A groupoid approach to discrete inverse semigroup algebras, Adv. Math. \textbf{223} (2010), 689--727.

\bibitem{St18} B.~Steinberg, Chain conditions on \'etale groupoid algebras with applications to Leavitt path algebras and inverse semigroup algebras, J. Aust. Math. Soc. \textbf{104} (2018), 403--411.

\bibitem{WZ19} L.~van Wyk, M.~Ziembowski, Lie solvability in matrix algebras, Linear Multilinear Algebra \textbf{67} (2019), 777--798.

\bibitem{ZS82} A.~E.~Zalesski\u{\i}, M.~B.~Smirnov, Associative rings satisfying the identity of Lie solvability, Vestsi Akad. Navuk BSSR Ser. Fiz.-Mat. Navuk 1982, no.~2, 15--20 (in Russian).

\end{thebibliography}
\end{document}